\input ./style/arxiv-general.cfg
\documentclass[MSNbibl,number,citesort,seceqn,dvips]{arxbj}
\makeatletter
   \@ifpackageloaded{graphicx}{}{\usepackage{graphicx}}
\makeatother
\usepackage{upgreek}

\volume{22}
\issue{3}
\pubyear{2016}
\firstpage{1331}
\lastpage{1363}
\doi{10.3150/14-BEJ694}
\docsubty{FLA}

\makeatletter
\newtheorem{theorem}{Theorem}[section]
\newtheorem{lemma}{Lemma}[section]
\newremark{remark}{Remark}[section]
\makeatother

\begin{document}
\begin{frontmatter}

\title{Goodness of fit tests in terms of local levels
with special emphasis on higher criticism~tests}
\runtitle{Goodness of fit tests in terms of local levels}

\begin{aug}
\author[A]{\inits{V.}\fnms{Veronika}~\snm{Gontscharuk}\corref{}\thanksref{A,B,e1}\ead[label=e1,mark]{veronika.gontscharuk@ddz.uni-duesseldorf.de}},
\author[A]{\inits{S.}\fnms{Sandra}~\snm{Landwehr}\thanksref{A,B,e2}\ead[label=e2,mark]{sandra.landwehr@ddz.uni-duesseldorf.de}}
\and\\
\author[B]{\inits{H.}\fnms{Helmut}~\snm{Finner}\thanksref{B,e3}\ead[label=e3,mark]{finner@ddz.uni-duesseldorf.de}\ead[label=u1,url]{www.foo.com}}
\address[A]{Department of Statistics in Medicine, Faculty of Medicine,
Heinrich-Heine-University, D\"usseldorf, Germany.\\ \printead{e1,e2}}
\address[B]{Institute for Biometrics and Epidemiology,
German Diabetes Center at the Heinrich-Heine-University, D\"usseldorf,
Germany. \printead{e3}}
\end{aug}
%
%
\received{\smonth{2} \syear{2014}}
%
\revised{\smonth{10} \syear{2014}}

%
\begin{abstract}
Instead of defining goodness of fit (GOF) tests in terms of their test
statistics, we
present an alternative method by introducing the concept of local
levels, which indicate high or low local sensitivity of a test.
Local levels can act as a starting point for the construction of new
GOF tests.
We study the behavior of local levels when applied to some well-known
GOF tests such as
Kolmogorov--Smirnov (KS) tests, higher criticism (HC) tests and tests
based on phi-divergences.
The main focus is on a rigorous characterization of the
asymptotic behavior of local levels of the original HC tests
which leads to several further asymptotic results for local levels of
other GOF tests
including GOF tests with equal local levels.
While local levels of KS tests, which are related to the central range,
are asymptotically strictly larger than zero, all local levels of HC tests
converge to zero as the sample size increases. Consequently, there
exists no asymptotic level $\alpha$ GOF test
such that all local levels are asymptotically bounded away from zero.
Finally, by means of numerical computations we compare classical KS and
HC tests to a GOF test with equal local levels.
\end{abstract}

%
\end{frontmatter}

\section{Introduction}\label{secintro}
Let $ X_1, \ldots, X_n $ be real-valued independently identically
distributed (i.i.d.) random
variables with continuous cumulative distribution function (c.d.f.) $F$.
We are interested in testing the null hypothesis
\begin{equation}
\label{nullhypothesis} 
H_0^{ \leq} \dvt F(x) \leq
F_0(x) \quad\mbox{or}\quad H_0^= \dvt F(x) =
F_0(x) \qquad \mbox{for all } x\in\mathbb{R},
\end{equation}
for a prespecified continuous c.d.f. $F_0$.
Since $F_0(X_i)$, $i=1,\ldots,n$, are i.i.d. uniformly distributed on
$[0,1]$ if $F = F_0$, we restrict our attention to the case where
\[
F_0 (x) = x \qquad \mbox{for all } x \in[0,1].
\]
Consequently, we assume that $X_i$, $i=1,\ldots,n$, take values in $[0,1]$.
We focus on the following class of goodness of fit (GOF) tests in terms
of order statistics $X_{1:n},\ldots,X_{n:n}$ related to the underlying
sample $X_1, \ldots, X_n $.
For testing $H_0^{ \leq}$ we consider a one-sided test $\varphi
\dvtx [0,1]\to\{0,1\}$ based on critical values $0 \leq c_{1,n} < \cdots<
c_{n,n}<1 $ such that
\begin{equation}
\label{defGOF} \varphi= 1 \quad \mbox{iff}\quad X_{i:n} \leq
c_{i,n} \qquad \mbox{for at least one } i=1,\ldots,n.
\end{equation}
A two-sided test $\tilde{\varphi} \dvtx  [0,1]\to\{0,1\}$ for testing
$H_0^=$ is defined by
\begin{equation}
\label{defGOF2} \tilde{\varphi} = 1\quad \mbox{iff}\quad X_{i:n}\leq
c_{i,n}\quad \mbox{or}\quad X_{i:n}\geq\tilde{c}_{i,n}
\qquad \mbox{for at least one } i=1,\ldots,n,
\end{equation}
where $0 \leq c_{1,n} < \cdots< c_{n,n}<1 $ and $0<\tilde{c}_{1,n} <
\cdots< \tilde{c}_{n,n} \leq1 $ are the corresponding critical values
fulfilling $c_{i,n}<\tilde{c}_{i,n}$, $i=1,\ldots,n$.
Thereby, $H_0^{\leq}$ is rejected if $\varphi=1$, while $H_0^{=}$ is
rejected if $\tilde{\varphi}=1$.
The global level of the test $\varphi$ and/or $\tilde{\varphi}$ is
given by $\mathbb{E}_0(\varphi) \equiv\mathbb{P}(\varphi=1 |
H_0^=)$ and/or $\mathbb{E}
_0(\tilde{\varphi}) \equiv\mathbb{P}(\tilde{\varphi}=1 |
H_0^=)$,
respectively.

We restrict attention to non-parametric tests only.
Among the most famous non-parametric GOF tests we find the
Kolmogorov--Smirnov (KS), Anderson--Darling (AD),
Cram\'er--von Mises and Berk--Jones (BJ) tests, where KS and BJ tests and
the supremum form of AD tests can be rewritten in the
form (\ref{defGOF}) and/or (\ref{defGOF2}).
In addition, recently proposed GOF tests based on the so-called
phi-divergences introduced in \cite{jagerwellner} are
non-parametric tests and can also be represented in the desired form.

In Section~\ref{sec11}, we briefly discuss the union-intersection principle
in relation to GOF tests and local levels.
Section~\ref{sec12} is concerned with the behavior of local levels of the
Kolmogorov--Smirnov test. Some further brief remarks concerning GOF
tests in terms of local levels
are given in Section~\ref{sec13}.
In Section~\ref{sec14}, we discuss the idea of GOF tests with equal
local levels, related ideas
and relations to recent work.
In Section~\ref{sec15}, we switch to higher criticism (HC) tests and
some further tests based on phi-divergences
and provide some figures which roughly illustrate the behavior
of the local levels of these tests.
An outline of the remaining part of the paper
the focus of which is on the asymptotics of local levels of the
original HC statistic
is given in Section~\ref{sec16}.

\subsection{The union-intersection principle and local levels} \label{sec11}

In multiple hypotheses testing local levels appear in a natural way, especially
in the case of multiple test procedures based on the union-intersection principle.
Such tests accept the global null hypothesis, that is, the intersection
of a
suitable set of elementary hypotheses $H_i$, if and only if all
elementary hypotheses are accepted.
Roughly speaking, a local level $\alpha_i $ for $H_i$
denotes the probability to reject $H_i$ if it is true.
Local levels tell us which amount of the overall level $\alpha$ is
attributed to each $H_i$.
Often multiple test procedures based on the union-intersection
principle have equal local levels.
Prominent examples are the classical Bonferroni test, Tukey's multiple
range test for pairwise comparisons, Dunnett's test for multiple
comparisons with a control or Scheffe's multiple contrast test.
A further general example is the minimum $p$-value test which
corresponds to the \textit{minimum level attained} test studied in \cite{berk78}.
The weighted Bonferroni test may serve as an example with different
local levels.

GOF tests of the form (\ref{defGOF}) and (\ref{defGOF2}) are related
to the union-intersection principle in the following way.
Let $U_1,\ldots,U_n$ be i.i.d. uniformly distributed random variables
and $U_{1:n},\ldots,U_{n:n}$ be the corresponding order statistics.
For each $i=1,\ldots,n$ consider null hypotheses $H_i^{ \leq}$ and
$H_i^{ =}$ on the distribution of a single order statistic $X_{i:n}$
such that $H_i^{ \leq}$ is true if $\mathbb{P}(X_{i:n}\leq x) \leq
\mathbb{P}
(U_{i:n}\leq x)$ for $x\in[0,1]$, that is, if $X_{i:n}$ is
stochastically larger than or equal to $U_{i:n}$ and $H_i^{ =}$ is true
if $X_{i:n}$ is equal to $U_{i:n}$
in distribution.
Define tests for $H_i^\leq$ by
\[
\varphi_i = 1\quad \mbox{iff}\quad X_{i:n} <
c_{i,n}
\]
and tests for $H_i^=$ by
\[
\tilde{\varphi}_i = 1 \quad \mbox{iff} \quad X_{i:n}\leq
c_{i,n} \quad\mbox{or}\quad X_{i:n}\geq\tilde{c}_{i,n}.
\]
Then $ H_0^{ \leq} \subseteq\bigcap_{i=1}^n H_i^{ \leq}$, $H_0^=
\subseteq\bigcap_{i=1}^n H_i^=$, $\{ \varphi= 1 \} = \bigcup_{i=1}^n \{
\varphi_i = 1 \}$ and $\{ \tilde{\varphi} = 1 \} = \bigcup_{i=1}^n
\{ \tilde{\varphi}_i = 1\}$ so that the GOF tests $\varphi$ and
$\tilde{\varphi}$ can be seen as union intersection tests.
We define \textit{local levels} of a GOF test by
\begin{equation}
\label{deflocallevels1}
\alpha_{i,n} = \mathbb{P} \bigl( \varphi_i =
1 | H_0^= \bigr) = \mathbb {P}(U_{i:n} \leq
c_{i,n})
\end{equation}
in the one-sided case and
\begin{equation}
\label{deflocallevels2}
\alpha^{=}_{i,n} = \mathbb{P} \bigl( \tilde{
\varphi}_i = 1 | H_0^= \bigr) = \mathbb{P}(
U_{i:n}\leq c_{i,n} ) + \mathbb{P}( U_{i:n}\geq
\tilde{c}_{i,n} )
\end{equation}
in the two-sided case.
Noting that $U_{i:n}$ is beta-distributed with parameters $i$ and
$n-i+1$ and denoting the related c.d.f. by $F_{i,n-i+1}$, we get $\mathbb{P}(
U_{i:n}\leq x ) = F_{i,n-i+1}(x)$.

Local levels can be viewed as an interesting characteristic of a GOF
test and may be interpreted as weights for
testing the family of null hypotheses $H_i^{\leq} $ or $H_i^=$,
$i=1, \ldots, n$.
The larger a local level $\alpha_{i,n}$ or $\alpha_{i,n}^= $, the
higher the chance to reject
the null hypothesis corresponding to the $i$th smallest order statistic
$X_{i:n}$ at least under the null hypothesis.
In other words, local levels can be regarded as a tool to signify areas
of high/low sensitivity of a test.
For example, if deviations from $H_0^\leq$ and/or $H_0^=$ are expected
in the tails, one would prefer a GOF test with larger local levels for
indices $i$ close to $1$ and/or close to $n$.
However, we have to take into account that order statistics are
dependent, see, for example, \cite{davidnagaraja} and \cite{shorackwellner}.
This may influence the probability to reject null hypotheses
corresponding to a set of $i$th
order statistics with indices $i$ in several ranges.

\subsection{Local levels of the Kolmogorov--Smirnov test} \label{sec12}
One of the most widely-used GOF tests, which can be written in terms of
(\ref{defGOF}) and/or (\ref{defGOF2}), is the well-known
Kolmogorov--Smirnov (KS) test.
We consider a one-sided asymptotic level $\alpha$  KS test, which
rejects $H_0^{\leq}$ if the KS test statistic
\[
\mathrm{KS}^+ = \max_{1\leq i\leq n} \sqrt{n} (i/n - X_{i:n})
\]
is larger than the asymptotic critical value $c_\alpha=\sqrt{-\log
(\alpha)/2}$ with $\alpha\in(0,1)$.
It holds $\lim_{n\to\infty}\mathbb{P}(\mbox{KS}^+ > c_\alpha|
H_0^=) = \alpha$, cf., for example,
\cite{shorackwellner}, page~11.
Even for $n\geq40$, the probability $\mathbb{P}(\mathrm{KS}^+ >
c_\alpha| H_0^=)$
is approximately $\alpha$.
The one-sided KS test can be represented in the form~(\ref{defGOF})
with critical values $c_{i,n}^{\mathrm{KS}} = \max(0,i/n-c_\alpha/ \sqrt
{n})$, $i=1,\ldots,n$.
In accordance with~(\ref{deflocallevels1}), the corresponding local
levels are given by $ \alpha_{i,n}^{\mathrm{KS}} = F_{i,n-i+1}(
c_{i,n}^{\mathrm{KS}})$, $i=1,\ldots,n$.
Note that $\alpha_{i,n}^{\mathrm{KS}}=0$ for $i\leq c_\alpha\sqrt{n}$.
For a finite $n$, the remaining $\alpha_{i,n}^{\mathrm{KS}}$ can be
calculated numerically.
Moreover, using the normal approximation, we get for $i\equiv i_n$
satisfying $i_n/n\to\zeta\in(0,1)$ that
\[
\lim_{n\to\infty} \alpha_{i_n,n}^{\mathrm{KS}} = 1-\Phi\bigl(
\sqrt{ -\log (\alpha ) /\bigl( 2 \zeta(1-\zeta)\bigr)} \bigr),
\]
where $\Phi(\cdot)$ is the standard normal cumulative distribution function.
The largest asymptotic local level is attained at $\zeta=1/2$ and
equals $1-\Phi(\sqrt{-2 \log\alpha})$.
\begin{figure}

\includegraphics{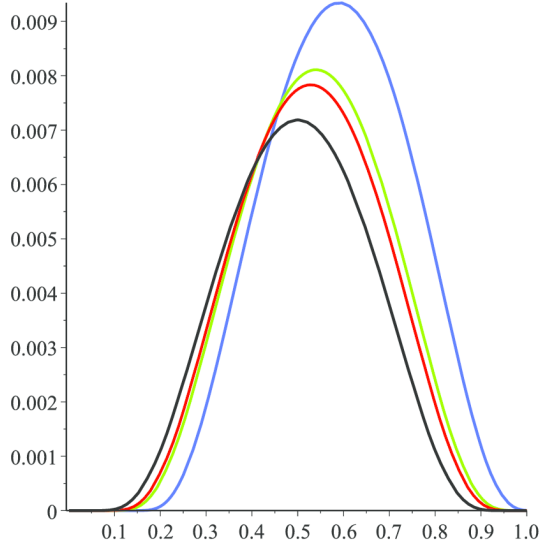}

\caption{Local levels $\alpha_{i,n}^{\mathrm{KS}}$ as a function of $i/n$
for one-sided KS tests with $\alpha=0.05$ and
$n=100, 500,1000$ together with the corresponding asymptotic local
levels (from top to bottom in
$i/n=0.8$).} \label{figlocallevelsKS}
\end{figure}
Figure~\ref{figlocallevelsKS} shows asymptotic and exactly
calculated local levels $\alpha^{\mathrm{KS}}_{i,n}$ as a function of $i/n$
for various $n$-values.
For $i_n/n$ in a central range of $[0,1]$, the limiting local levels
are bounded away from zero,
whereas for $i_n/n\to\zeta\in\{0,1\}$ we get $\lim_{n\to\infty
}\alpha_{i_n,n}^{\mathrm{KS}}=0$.
This coincides with the well-known fact that KS tests have higher power
for alternatives that differ from the null distribution in the central
range and low power against alternative distributions which mainly deviate
from the null in the tails. Alternatives of this kind, however, are
common in many applications,
for example, in genome-wide association studies, in which we face a
very large number of hypotheses to
test with only a small number of them being non-null.
For more practical applications see, for example, \cite{cai07,halljin08} and \cite{halljin10}.

Various modifications of the KS test have been proposed in the past.
For example, R{\'e}v{\'e}sz \cite{Revesz} constructed a test based on a
statistic which combines the advantages of the classical and normalized
KS statistics with regard to their sensitivity ranges. Mason and
Schuenemeyer \cite{mason} introduced a
modified KS test by combining the classical KS with R{\'e}nyi-type
statistics and investigated the finite sample and asymptotic
distribution of this modification. Test statistics that are determined
by order statistics, in particular
tail order statistics, are studied by Lockhart in \cite{lockhart} with
respect to asymptotic relative efficiency against
a certain class of alternatives. Bahadur efficiencies for a lot of
non-parametric GOF tests are extensively studied by Nikitin in \cite{nikitin}.
More recently, Jager and Wellner \cite{jagerwellner} proposed GOF
tests based on phi-divergences. Their supremum- and integral-type
statistics cover various forms of Anderson--Darling and Berk--Jones statistics.

\subsection{GOF tests in terms of local levels} \label{sec13}
For many (non-parametric) GOF tests, there is a class of alternatives
against which this test is the most powerful.
Hence, if we have some information about the range, where the
alternative distribution mainly deviates from the null distribution, it
seems worthwhile to apply such an appropriately tailored GOF test.
However, from the viewpoint of test statistics it is difficult to
determine whether the corresponding GOF test is sensitive for a
predefined range of deviations.
Fortunately, the construction of tailored GOF tests is much easier by
means of local levels.
Thereby, the aim is to construct a GOF test
with larger local levels in the crucial area.
For example, assuming sparse signals, a GOF test with larger local
levels for indices close to $1$ and smaller local levels for the
remaining indices seems to be a reasonable choice.

In general, for a given suitable set of local levels $\alpha_{i,n}$,
$i=1,\ldots,n$, we are able to construct a corresponding GOF test of
the form (\ref{defGOF}) and/or (\ref{defGOF2}).
In the one-sided case the related critical values are given by $c_{i,n}
= F_{i,n-i+1}^{-1}(\alpha_{i,n})$, $i=1,\ldots,n$, where $F_{i,n-i+1}^{-1}$
denotes the inverse function of the c.d.f. $F_{i,n-i+1}$.
For a two-sided GOF test $\tilde{\varphi}$ we have to decide how to
split $\alpha_{i,n}^=$ into two non-negative terms $\alpha_{i,n}^{(1)}$
and $\alpha_{i,n}^{(2)}$ such that $\alpha_{i,n}^{(1)} + \alpha
_{i,n}^{(2)} = \alpha_{i,n}^=$ and $\mathbb{P}( U_{i:n} \leq c_{i,n}
)=\alpha
_{i,n}^{(1)}$, $\mathbb{P}(U_{i:n} \geq\tilde{c}_{i,n} )=\alpha
_{i,n}^{(2)}$.
One possibility may be $\alpha_{i,n}^{(1)} = \alpha_{i,n}^{(2)} =
\alpha
_{i,n}^= /2$, which leads to $c_{i,n} = F_{i,n-i+1}^{-1}(\alpha
_{i,n}^=/2)$ and $\tilde{c}_{i,n} = F_{i,n-i+1}^{-1}(1-\alpha_{i,n}^=/2)$.
The latter can be calculated at least numerically.

\subsection{GOF tests with equal local levels} \label{sec14}
If we do not have any idea on alternatives, it seems natural to choose
a GOF test with equal local levels, that is,
\[
\alpha_{1,n} = \cdots= \alpha_{n,n} = \alpha_n^{\mathrm{loc}}
\quad\mbox{and/or}\quad \alpha_{1,n}^= = \cdots= \alpha_{n,n}^= =
\alpha _n^{\mathrm{loc}}
\]
for some suitable $\alpha_n^{\mathrm{loc}}\in(0,1)$.
The idea behind this proposal is similar to the idea behind the KS
test, where the distance between the empirical c.d.f. $\hat{F}_n(x)$ and
the underlying c.d.f. $F_0(x)$, that is, $\hat{F}_n(x)-F_0(x)$ for the
one-sided test case and $|\hat{F}_n(x)-F_0(x)|$ for the two-sided case,
is compared to the same critical value for each $x$.
That is, the KS test can be seen as a GOF test with equal distances for
all feasible $x$-values.
Considering other measures of the distance between the theoretical and
the corresponding empirical distributions, one may construct various
GOF tests with some quantities being equal.
For example, a family of GOF tests introduced in \cite{jagerwellner}
can be seen as tests with equal phi-divergences.
A prominent example here is the Berk--Jones test which corresponds to
equal Kullback-Leibler divergences.
Altogether, the idea of considering equal quantities such as equal
distances, critical values, test statistics and also local levels, is a
natural approach when constructing GOF tests.

GOF tests with local levels equal to some $\alpha_n^{\mathrm{loc}}\in(0,1)$
are given as follows.
The one-sided version of the test $\varphi(\alpha_n^{\mathrm{loc}})$ (say)
is defined by (\ref{defGOF}) with $c_{i,n}=F_{i,n-i+1}^{-1}(\alpha
_n^{\mathrm{loc}})$, $i=1,\ldots,n$.
The two-sided test $\tilde{\varphi}(\alpha_n^{\mathrm{loc}})$
is given
by (\ref{defGOF2}) with $c_{i,n}=F_{i,n-i+1}^{-1}(\alpha_n^{\mathrm
{loc}}/2)$ and $\tilde{c}_{i,n}=1-F_{i,n-i+1}^{-1}(\alpha_n^{\mathrm
{loc}}/2)$, $i=1,\ldots,n$.
In order to get a level $\alpha$ test we have to choose $\alpha
_n^{\mathrm{loc}}$ such that $\mathbb{E}_0(\varphi(\alpha
_n^{\mathrm{loc}}) ) = \alpha$ and/or $\mathbb{E}
_0(\tilde{\varphi}(\alpha_n^{\mathrm{loc}}) ) = \alpha$.
Unfortunately, it seems there does not exist any analytically
manageable formula for $\alpha_n^{\mathrm{loc}}$ as a function of $n$ and
$\alpha$ so that $\alpha_n^{\mathrm{loc}}$ has to be calculated numerically.
Nevertheless, we are able to provide some bounds for $\alpha
_n^{\mathrm{loc}}$.
For example, the Bonferroni inequality implies
\[
\alpha/ n < \alpha_n^{\mathrm{loc}} < \alpha,\qquad  n\in\mathbb{N}.
\]
Moreover, it can easily be seen that $\alpha_n^{\mathrm{loc}}$ lies between
the smallest and largest local levels for any (exact) level $\alpha$
GOF test of type (\ref{defGOF}) and (\ref{defGOF2}), respectively.
Thus, knowledge of local levels corresponding to suitable GOF tests
leads at least to upper and lower bounds for $\alpha_n^{\mathrm{loc}}$.
For example, by means of the asymptotic KS local levels, we get for the
one-sided case
\[
0 < \alpha_n^{\mathrm{loc}} \leq\Phi\bigl(\sqrt{-2\log(\alpha)}
\bigr) + \mathrm{o}(1),\qquad n\in\mathbb{N},
\]
which is, unfortunately, a very wide range.
Thus, we have to study local levels related to other level $\alpha$
GOF tests.

Once we have $\alpha_n^{\mathrm{loc}}$, one may redefine
the corresponding GOF tests with equal local levels
as minimum $p$-value (minP) tests based on the one-sided $p$-values
$p_{i,n} = F_{i,n} ( X_{i:n} )$, $ i=1,\ldots,n$.
Setting $ M_{n}^+ = \min_{1\leq i\leq n} p_{i,n} $
and $ M_{n} = \min_{1 \leq i \leq n} \{ p_{i,n}, 1- p_{i,n} \} $,
we get $ \varphi(\alpha_n^{\mathrm{loc}}) = 1 $ iff $ M_{n}^+ \leq
\alpha
_n^{\mathrm{loc}}$
and $ \tilde{\varphi}(\alpha_n^{\mathrm{loc}}) = 1$ iff
$M_{n} \leq
\alpha_n^{\mathrm{loc}}/2$.

The minP statistics $ M_{n}^+$ and $ M_{n} $
were already introduced
by Berk and Jones in 1979 (cf. \cite{berk79})
and they referred to these statistics
as minimum level attained statistics.
Implicitly, Berk and Jones
were the first proposing the construction of
equal local level GOF tests (even though they did not use the term
local levels).
Among others, they
extensively studied $M_{n}^+$ and $ M_{n} $
with respect to optimality and Bahadur efficiency, see also \cite{berk78}.
A further representation of GOF tests with equal local levels was
provided in the unpublished manuscript \cite{buja06} in 2006.
Moreover, such tests were recently provided by several authors.
At the 7th International Conference on Multiple Comparison Procedures
(MCP) 2011 we introduced the concept of local levels and proposed GOF
tests with equal local levels
as an improvement of the higher criticism (HC) tests.
At the MCP 2013 we presented asymptotic as well as finite properties of
GOF tests
with equal local levels, cf. \cite{GLF2014}.
In contrast to the formulation via local levels, the GOF test in \cite
{buja06} is formulated in terms of bounding functions.
This representation of the test is elaborated on in \cite{aldor13}.
What is more, the same test is provided in \cite{mary14} in the HC framework.
Finally, the test is also considered in the preprints \cite{kaplan},
\cite{kaplannew} and \cite{moscovich}.

\subsection{Higher criticism and phi-divergence} \label{sec15}
In connection with high dimensional data and associated multiple
testing issues, the so-called higher criticism (HC) tests generated
considerable interest during the last decade, cf. for example, \cite{donohojin04,donohojin08,donohojin09} and
\cite{halljin10}.
For example, Donoho and Jin proposed the use of HC tests when testing
the global null hypothesis in high-dimensional
models with sparse signals against some specific alternatives.
Studying the HC test statistic they showed in \cite{donohojin04}
that, asymptotically, HC related tests are successful throughout the
same region of amplitude sparsity where the corresponding oracle
likelihood ratio test would succeed.
This means that a further specification of an alternative is not necessary.
Note that HC tests can also be seen as GOF tests of the type (\ref
{defGOF}) and/or (\ref{defGOF2}).
What is more, it appears that studying  local levels corresponding to
HC tests is essential in order to construct new GOF tests, which have a
high power against alternative distributions that mainly deviate from
the null distribution in the considered range.

\begin{figure}

\includegraphics{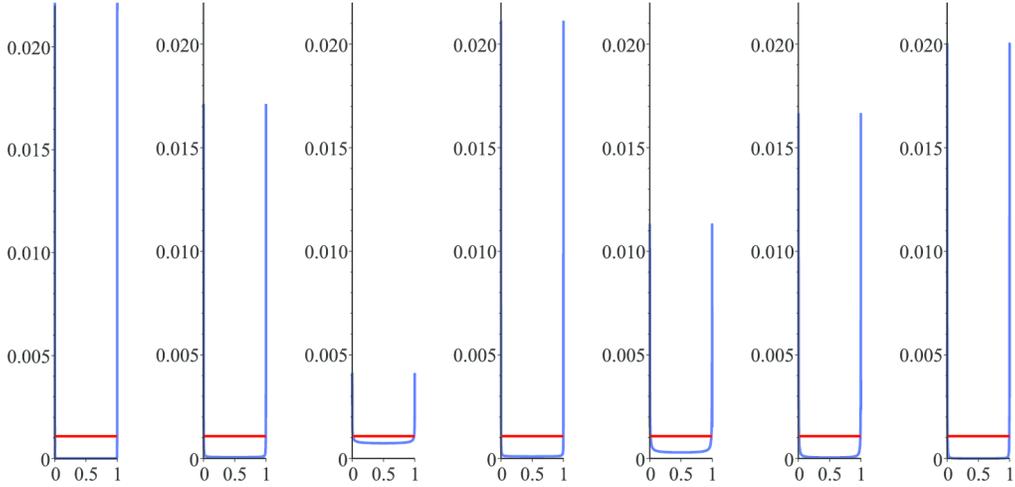}

\caption{Two-sided local levels of level $\alpha$ GOF tests based on
phi-divergences with $s=2,1.5,1,0.5,0, -0.5,-1$ (from left to right)
together with $\alpha_n^{\mathrm{loc}}=0.001075$ related to the two-sided
level $\alpha$ test $\tilde{\varphi}(\alpha_n^{\mathrm{loc}})$ for
$\alpha=0.05$ and $n=1000$.} \label{figloclevphidiv}
\end{figure}

Alternatively, instead of HC tests one may consider other GOF tests
which are based on the phi-divergences introduced in \cite{jagerwellner}.
Thereby, the family of these tests is parametri\-zed by $s\in[-1,2]$ so
that the HC test corresponds to $s=2$, the Berk--Jones test to $s=1$,
the reversed Berk--Jones test to $s=0$ and the studentized version of
the HC test to $s=-1$.
As suggested by a referee, we compare local levels for some selected $s$-values.
Figure~\ref{figloclevphidiv} shows two-sided local levels of the
exact level $\alpha$ tests based on the phi-divergences for $\alpha
=0.05$, $n=1000$ and $s=2,1.5,1,0.5,0,-0.5,-1$.
What these local levels have in common is that they are large in the
tails and small and approximately equal in the central range.
However, the range of the local level values is largest for the HC test
and smallest for the Berk--Jones test.
Therefore, it looks that the Berk--Jones test leads to the narrowest
bounds for $\alpha_n^{\mathrm{loc}}$ while the HC test to the widest ones.
Due to the fact that under the null hypothesis statistics based on
phi-divergences have the same asymptotic behavior in a specific range
relevant for the asymptotics, any of these tests will lead to the same
asymptotic results for most local levels.
Therefore, it does not matter which test we consider.
Since the tests with $s=2$ (HC tests) and $s=-1$ (studentized HC tests)
have the simplest representation of the form (\ref{defGOF}) and/or
(\ref{defGOF2}), we prefer to restrict attention to the original HC
test, which has received a lot of attention during the past decade.
\subsection{Outlook of the remaining part of the paper} \label{sec16}
In this paper, we calculate local levels of asymptotic level $\alpha$
HC tests and show that these local levels converge to zero as $n\to
\infty$, which differs drastically from the KS case, cf. Figure~\ref{figlocallevelsKS}.
This implies for local levels of any asymptotic level $\alpha$ GOF test
of the form (\ref{defGOF}) and/or (\ref{defGOF2}) that
\[
\lim_{n\to\infty}\min_{1\leq i \leq n} \alpha_{i,n}
=0,
\]
that is, there are no level $\alpha$ tests for which the local levels
are all asymptotically bounded away from zero.
Finally, by a careful study of asymptotic HC local levels we get for
$\varphi\equiv\varphi(\alpha_n^{\mathrm{loc}})$ and/or $\varphi
\equiv
\tilde{\varphi}(\alpha_n^{\mathrm{loc}})$ that
\[
\lim_{n\to\infty}\mathbb{E}_0(\varphi) = \alpha \quad \mbox{iff}\quad \lim_{n\to\infty}\alpha_n^{\mathrm{loc}} \cdot
\frac{2 \log(\log
(n)) \log(n) }{-\log(1-\alpha)}=1.
\]
This result seems to be the most precise result concerning the asymptotics
of the one- and two-sided GOF tests with equal local levels.
In general, there are only few other works, in which asymptotics is
investigated, cf.  \cite{kaplan,kaplannew} and \cite{moscovich}.
Due to a long revision process, some highlights of this paper have been
already summarized in \cite{GLF2014}, where the focus lies on the
sensitivity range of the HC tests statistic, extremely slow HC
asymptotics, relations to the Ornstein--Uhlenbeck process, and power
comparisons of the test with equal local levels and the original HC test.
The remaining part of the paper is organized as follows.
In Section~\ref{secHC}, we study local levels of the HC test.
We further derive the critical value and rejection curves corresponding
to asymptotic level~$\alpha$ HC tests and provide a result on the
asymptotic behavior of the HC critical values.
As zones of normal and Poisson convergence play a crucial role in the
derivation of asymptotic results,
we provide some basic results on these approximations for HC local
levels in Section~\ref{secapproximations}.
Section~\ref{secexpressions} contains explicit asymptotic expressions
of the local levels $\alpha_{i,n}$ of
the one-sided HC test. They are derived for various growth rates of $i$
utilizing the approximation results
from Section~\ref{secapproximations}. In Section~\ref{secquotient},
we investigate the asymptotic
monotonicity of the local levels of one-sided HC tests and provide some
results concerning the asymptotic behavior of local levels related to
general level $\alpha$ GOF tests and tests with equal local levels.
In Section~\ref{secGOFcomparison}, we compare classical KS and HC
tests to GOF tests with equal local levels by means of numerical computations.
Future investigations and open questions are discussed in Section~\ref{outlook}.
All proofs mostly of technical nature are deferred to Appendices~\ref{AppA}, \ref{AppB} and \ref{AppC}.
\section{Higher criticism tests and local levels} \label{secHC}
First, we introduce the version of the higher criticism GOF tests that
we are dealing with.
Let
\[
G_{i,n} ( u ) = \sqrt{ n } \frac{i/n - u}{\sqrt{u ( 1 - u ) } }\quad \mbox{and}\quad
\tilde{G}_{i,n} ( u ) = \sqrt{ n } \frac{ u -(i-1)/n }{\sqrt{ u
(1-u) } },\qquad u \in(0, 1).
\]
A class of one-sided and two-sided HC test statistics can be expressed as
\[
\mbox{HC}^+ = \max_{ 1 \leq i \leq n } G_{i,n} ( X_{i:n} )
\quad \mbox{and} \quad \mbox{HC}^{=} = \max_{ 1 \leq i \leq n } \bigl\{
G_{i,n} ( X_{i:n} ), \tilde {G}_{i,n} (
X_{i:n} ) \bigr\},
\]
respectively.
A one-sided HC test based on a critical value $d$ rejects $H_0^\leq$
iff $\mbox{HC}^+ > d$ and a two-sided HC test with the same critical value
rejects $H_0^{=}$ iff $ \mbox{HC}^{=} > d$.
%
In accordance with (\ref{deflocallevels1}), local levels of a
one-sided HC test based on a critical value $d>0$ are given by
\begin{equation}
\label{deflocallevelsonesided}
\alpha_{i,n} = \mathbb{P} \bigl( G_{i,n} (
U_{i:n} ) > d \bigr),\qquad  i = 1,\ldots,n.
\end{equation}
Analogously, local levels of the corresponding two-sided HC test are
given by
\begin{equation}
\label{deflocallevelstwosided}
\alpha^{=}_{i,n} = \mathbb{P} \bigl( \bigl\{
G_{i,n} ( U_{i:n} ) > d \bigr\} \cup\bigl\{
\tilde{G}_{i,n} ( U_{i:n} ) > d \bigr\} \bigr), \qquad i=1,\ldots,n.
\end{equation}
%
Setting $ u_n = \log(\log(n)) $ and
\begin{equation}
\label{dn} d_n(t) = \bigl( t + 2 u_n + \bigl( \log(
u_n) - \log(\uppi) \bigr) / 2 \bigr) / \sqrt { 2 u_n },
\end{equation}
the asymptotic distributions of $\mbox{HC}^+$ and $\mbox{HC}^=$ are
given by
\begin{equation}
\label{wholeHC}
\lim_{n \to\infty} \mathbb{P}\bigl( \mbox{HC}^+ \leq
d_n(t) | H_0^= \bigr) = \exp\bigl( - \exp( - t ) \bigr)
\end{equation}
and
\begin{equation}
\label{wholeHCbetrag}
\lim_{n \to\infty} \mathbb{P}\bigl(
\mbox{HC}^{=} \leq d_n(t) | H_0^= \bigr) = \exp
\bigl( - 2 \exp( - t ) \bigr),
\end{equation}
respectively, cf. \cite{eicker79} and \cite{jaeschke79}.
For $t=t_\alpha$ or $t=t_\alpha^=$ with
\[
t_{\alpha}=-\log\bigl(-\log(1-\alpha)\bigr) \quad\mbox{and}\quad
t_{\alpha}^{=} = -\log\bigl( -\log(1-\alpha)/2 \bigr)
\]
we get a one-sided or two-sided asymptotic level $\alpha$ HC test,
respectively.
Note that
\[
\bigl\{ G_{i,n} ( U_{i:n} ) > d \bigr\} = \bigl\{
U_{i:n} < h_{i,n}(d) \bigr\} \quad\mbox{and}\quad \bigl\{
\tilde{G}_{i,n} ( U_{i:n} ) > d \bigr\} = \bigl\{
U_{i:n} > \tilde{h}_{i,n}(d) \bigr\},
\]
where
\begin{equation}
\label{hn}
h_{i,n}(d) = {\frac{d^2 + 2 i - d \sqrt{d^2 + 4 i - 4 i^2 /n }}{ 2 (
d^2+n) }}
\end{equation}
and
\begin{equation}
\label{h2n}
\tilde{h}_{i,n}(d) = {\frac{d^2 + 2( i-1) + d \sqrt{d^2 + 4 (i-1) - 4
(i-1)^2 /n }}{ 2 ( d^2+n) }}.
\end{equation}
Thereby, $\tilde{h}_{i,n}(d) > h_{i,n}(d)$, $i=1,\ldots,n$, for $d$
large enough.
Below, let $d\geq1$, which guarantees $\tilde{h}_{i,n}(d) >
h_{i,n}(d)$ for all $i=1,\ldots,n$ and
$\{ U_{i:n} < h_{i,n}(d) \} \cap\{ U_{i:n} > \tilde{h}_{i,n}(d) \} =
\varnothing$.
Thus, local levels can be expressed as
\[
\alpha_{i,n} = \mathbb{P}\bigl( U_{i:n} < h_{i,n}
\bigl(d_n(t)\bigr) \bigr), \qquad i=1,\ldots,n,
\]
for one-sided HC tests and
\[
\alpha^{=}_{i,n} = \mathbb{P}\bigl( U_{i:n} <
h_{i,n}\bigl(d_n(t)\bigr) \bigr) + \mathbb {P} \bigl(
U_{i:n} > \tilde{h}_{i,n}\bigl(d_n(t)\bigr) \bigr),\qquad i=1,\ldots,n,
\]
for two-sided HC tests.
Assuming that $Z_n$ ($\tilde{Z}_n$) is a binomially distributed
random variable with parameters\vspace*{1pt} $n$ and $h_{i,n}(d_n(t))$ ($\tilde
{h}_{i,n}(d_n(t))$), that is,
$Z_n \sim\mathcal{B}(n,h_{i,n}(d_n(t)))$, and $\tilde{Z}_n \sim
\mathcal{B}(n,\tilde{h}_{i,n}(d_n(t)))$, we get
\begin{equation}
\label{binomialformula}
\alpha_{i,n} = \mathbb{P}(Z_n \geq i)\quad \mbox{and} \quad\alpha^{=}_{i,n} = \mathbb{P}(Z_n \geq i) +
\mathbb{P}(\tilde{Z}_n < i).
\end{equation}
Since $\tilde{h}_{i,n}(d)=1-h_{n,n-i+1}(d)$ for $i=1,\ldots,n$ and
$n\in
\mathbb{N}$, we obtain
\[
\alpha^{=}_{i,n} = \alpha_{i,n} +
\alpha_{n-i+1,n},\qquad i=1,\ldots,n.
\]
Hence, two-sided local levels are symmetric in the sense $\alpha
^{=}_{i,n} = \alpha^{=}_{n-i+1,n}$ for
$i=1,\ldots,n$, and can be easily calculated if one-sided local levels
are known.

\begin{figure}

\includegraphics{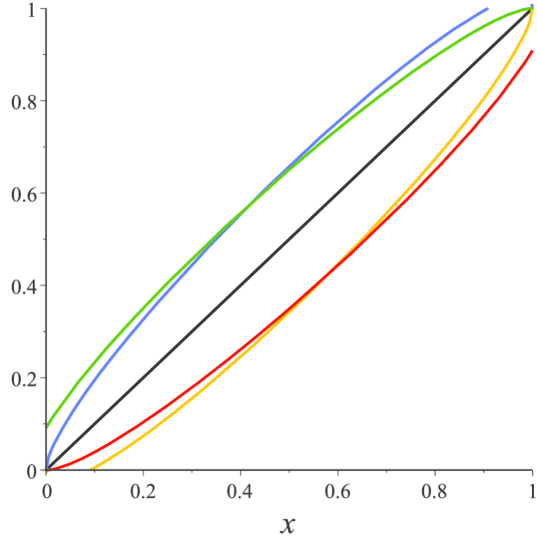}

\caption{For $n = 1000$ and $d = 10$ the critical value curve $\rho_n(x, d)$ and the corresponding
rejection curve $r_n(x, d)$ are given by the lowest and highest curves in $x = 0.8$,
respectively;
$\tilde{\rho}_n(x, d)$ and $\tilde{r}_n(x, d)$ are given by the highest and lowest curves in $x = 0.2$,
respectively. The straight line is the diagonal.}
\label{figcritvalandrejectcurves}
\end{figure}

Note that the considered HC tests (and a lot of multiple tests) can be
alternatively defined in
terms of a rejection curve, which is a general inverse of the
corresponding critical value curve, cf. \cite{FiDiRo2009}.
Critical value curves related to (\ref{hn}) and (\ref{h2n}) are given by
\begin{equation}
\label{rhon}
\rho_n(x,d) = \frac{d^2 + 2 x n - d \sqrt{d^2 + 4 x n - 4 x^2 n }}{ 2
( d^2+n) }
\end{equation}
and
\begin{equation}
\label{rho2n}
\tilde{\rho}_n(x,d) = \frac{d^2 + 2 (x n -1) + d \sqrt{d^2 + 4 (x n-1)
- 4 (x n -1)^2/n }}{ 2 ( d^2+n) },
\end{equation}
respectively, that is, $ h_{i,n}(d) = \rho_n(i/n,d) $ and $ \tilde
{h}_{i,n}(d) = \tilde{\rho}_n(i/n,d) $.
The corresponding rejection curves are given by
\[
r_n(x,d) = x + d \sqrt{ \frac{ x (1-x)}{n} } \quad\mbox{and}\quad
\tilde{r}_n(x,d) = x + \frac{1}{n} - d \sqrt{ \frac{ x (1-x)}{n}
},
\]
respectively.
It holds $\rho_n(x,d) = 1-\tilde{\rho}_n(1-x+1/n,d)$ and $r_n(x,d) =
1-\tilde{r}_n(1-x,d)+1/n$.
Figure~\ref{figcritvalandrejectcurves} shows critical value curves
$\rho_n$,
$\tilde{\rho}_n$ and the corresponding rejection curves $r_n$,
$\tilde
{r}_n$ for $n=1000$ and $d=10$.
For increasing $n$ and/or decreasing $d$, the corresponding curves tend
to the diagonal.

%

The following lemma shows the asymptotic behavior of the critical
values $h_{i_n,n}\equiv h_{i_n,n}(d_n(t))$ for different ranks $i$.

\begin{lemma}\label{entwicklungh}
Let $n\in\mathbb{N}$ and $i_n\in\mathbb{N}$ with $i_n\leq n$. It
holds:
\begin{longlist}[(iii)]
\item[(i)]  if $i_n=\mathrm{o}(u_n)$ as $n\to\infty$, then
\begin{eqnarray}
 \frac{ n h_{i_n,n} }{i_n} & = &
\frac{i_n}{2u_n} \biggl( 1 -
\frac{\log(u_n)+2t-\log(\uppi)}{2 u_n} - \frac{i_n}{ u_n} \biggr)
\nonumber
\\[-8pt]
\label{hseriesextremeleft}
\\[-8pt]
\nonumber
&& {}+ \mathrm{O} \biggl( \frac{i_n \log(u_n)^2 +i_n^3 }{u_n^3 } \biggr);
\end{eqnarray}
\item[(ii)] if $c_n\equiv u_n/i_n \to c$ for some $c>0$ we obtain
\begin{equation}
\label{hseriesintermediate}
\frac{ n h_{i_n,n} }{i_n} = \delta(c_n) \biggl( 1 -
\frac{\log(u_n)+2 t-\log(\uppi)}{2 i_n \sqrt
{c_n^2+2c_n}} \biggr) + \mathrm{O} \biggl( \frac{\log(u_n)^2}{u_n^2} \biggr),
\end{equation}
where $\delta(c_n)= 1+c_n-\sqrt{c_n^2+ 2 c_n}\in(0,1)$;
\item[(iii)]  if $i_n(1-i_n/n)/u_n\to\infty$, that is, $i_n/u_n\to\infty$ and
$(n-i_n)/u_n\to\infty$, then
\begin{eqnarray}
\frac{ n h_{i_n,n} }{i_n} & = &
1 - \sqrt{\frac{2 u_n}{i_n}
\biggl( 1- \frac{i_n}{n} \biggr)} - \frac{\log(u_n)
+2 t-\log(\uppi)}{2\sqrt{2 i_n u_n}} \sqrt{ 1-
\frac{i_n}{n} }
\nonumber
\\[-8pt]
\label{hseriescentralintermed}\\[-8pt]
\nonumber
&&{}+ \biggl(1-\frac{2 i_n}{n} \biggr) \frac{u_n}{i_n} + \mathrm{o} \biggl(
\frac{u_n}{i_n} + \frac{ 1}{\sqrt{i_n u_n} } \sqrt {1-\frac
{i_n}{n}} \biggr);
\end{eqnarray}
\item[(iv)]  if $c_n\equiv(n-i_n)/u_n \to c $ for some $c\geq0$ we obtain
\begin{eqnarray}
\frac{n h_{i_n,n}}{i_n} & = &
1 - \frac{u_n}{i_n} ( 1 +
\sqrt{1+2 c_n} ) -\frac{\log
(u_n)+2t-\log(\uppi)}{2 i_n}
\nonumber
\\[-8pt]
\label{hseriesextremeright}
\\[-8pt]
\nonumber
&&{}\times \biggl( 1+ \frac{1+c_n}{ \sqrt{1+2 c_n} } \biggr) + \mathrm{O} \biggl( \frac{\log(u_n)^2}{i_n u_n}
\biggr).
\end{eqnarray}
\end{longlist}
\end{lemma}

In the next sections, we provide asymptotics of local levels $\alpha
_{i,n}$ of HC tests for all growth rates of $i$.
To be precise, we are considering so-called extreme ranks $i$, where
$i$ or $n-i$ are fixed, and increasing ranks $i=i_n \to\infty$ as
$n\to
\infty$.
We split the latter into central ranks, which are such that $i_n/n \to
\zeta\in(0,1)$ as $n\to\infty$,
and intermediate ranks, which are such that $i_n/n \to\zeta\in\{0,1\}$.
For these concepts see, for example, \cite{leadbetter}.
\section{Normal and Poisson approximations for local levels}\label{secapproximations}
Due to representation (\ref{binomialformula}), we can approximate
local levels of a HC test by applying Poisson and/or normal
approximations for the binomial distribution.
Below, let $Y_n\sim\mathcal{P}(n h_{i_n,n})$ and $\tilde{Y}_n
\sim\mathcal{P}(n(1-h_{i_n,n}))$, where $\mathcal{P}(\lambda)$ denotes
the Poisson distribution with parameter $\lambda> 0$.
Thereby, $h_{i_n,n}=h_{i_n,n}(d_n(t))$ is given in Lemma~\ref{entwicklungh}.

The following theorem shows that for large values of $n$ local levels
$\alpha_{i_n,n}$ of HC tests
based on critical values $h_{i_n,n}$, can be calculated by means of
Poisson approximations
for a wide range of ranks $i_n$.

\begin{theorem}[(Poisson approximation of local levels)]\label{lemloclevPois}
Let $i_n\in\mathbb{N}$, $i_n\leq n$, be a sequence of non-decreasing numbers.
For $i_n$ such that $i_n=\mathrm{o}(u_n) $ we obtain
\begin{equation}
\label{eqloclevPois}
\alpha_{i_n,n} = \mathbb{P}(Y_n =
i_n) \bigl[1+\mathrm{o}(1) \bigr],
\end{equation}
for $i_n$ such that $u_n/i_n\to c $ for some $c>0$
\begin{equation}
\label{inCun} \alpha_{i_n,n} = 1/\bigl(\sqrt{c^2+2 c}-c\bigr)
\mathbb{P}(Y_n = i_n)\bigl[1+\mathrm{o}(1) \bigr]
\end{equation}
and for $i_n$ such that $i_n/u_n\to\infty$ and $i_n=\mathrm{o}(\sqrt{n/u_n} )$
\begin{equation}
\label{eqloclevPois2}
\alpha_{i_n,n} = \sqrt{i_n/(2 u_n)}
\mathbb{P}(Y_n = i_n) \bigl[1+\mathrm{o}(1) \bigr].
\end{equation}
Analogously, for $i_n$ with $n-i_n=\mathrm{o}(u_n)$, we get
\begin{equation}
\label{eqloclevPoisright}
\alpha_{i_n,n} = \mathbb{P}(\tilde{Y}_n =
n-i_n) \bigl[1+\mathrm{o}(1) \bigr],
\end{equation}
for $i_n$ fulfilling $(n-i_n)/u_n \to c $ for some $c>0$ we obtain
\begin{equation}
\label{ninCun}
\alpha_{i_n,n} = \bigl( 1 +c/(1+\sqrt{1+2 c}) \bigr)
\mathbb {P}(\tilde {Y}_n = n- i_n) \bigl[1 +\mathrm{o}(1) \bigr]
\end{equation}
and if $(n-i_n)/u_n\to\infty$ and $n-i_n=\mathrm{o}(\sqrt{n/u_n})$, then
\begin{equation}
\label{eqloclevPois22} \alpha_{i_n,n} = \sqrt{(n-i_n)/(2
u_n)} \mathbb{P}(\tilde{Y}_n = n- i_n)
\bigl[1+\mathrm{o}(1) \bigr].
\end{equation}
\end{theorem}

The following theorem shows that local levels of HC tests corresponding
to central ranks
and to intermediate ranks close to central ones can be calculated in
terms of the density of the standard normal  distribution $\phi$.

\begin{theorem}[(Normal approximation of local levels)]\label{zoneofconv}
Let $i_n\in\mathbb{N}$ be such that $i_n(1-i_n/n)/u_n^3 \to\infty$,
$ \sigma
_n =
\sqrt{n h_{i_n,n}(1-h_{i_n,n}) } $ with $h_{i_n,n}$ given in (\ref
{hseriescentralintermed}) and $x_n=(i_n-n h_{i_n,n})/\sigma_n$.
Then $x_n\to\infty$, $x_n^3/\sigma_n\to0$ as $n\to\infty$ and
\begin{equation}
\label{quotient2} \alpha_{i_n,n} = \phi(x_n)/ x_n
\bigl[ 1 + \mathrm{O} \bigl( x_n^3 / \sigma_n +
1/x_n^2 \bigr) \bigr].
\end{equation}
\end{theorem}

\begin{pf}
We can derive (\ref{quotient2}) by following the proof in \cite{smirnov},
where he considered the case $p_n \equiv p$.
Since\vspace*{1pt} $\alpha_{i_n,n}= \mathbb{P}(Z_n \geq i_n)$, where $Z_n \sim
\mathcal
{B}(n,h_{i_n,n})$,
it suffices to show $x_n^3/\sigma_n\to0$ if $i_n(1-i_n/n)/u_n^3 \to
\infty$.
This can easily be proved by applying (\ref{hseriescentralintermed}),
which implies
$\sigma_n = \sqrt{i_n(1-i_n/n)}[1+\mathrm{o}(1)] $ and $x_n=\sqrt{2 u_n}[1+\mathrm{o}(1)]$.
\end{pf}

\begin{remark}
Note that for $i_n$ satisfying $i_n(1-i_n/n)/u_n^3 \to\infty$ as
$n\to
\infty$ and
$i_n(1-i_n/n)=\mathrm{o}(\sqrt{n/u_n})$, Theorems \ref{lemloclevPois} and \ref
{zoneofconv}
provide two alternative approximations for local levels of HC tests.
\end{remark}

\section{Asymptotic expressions of local levels of HC tests} \label
{secexpressions}
\begin{figure}

\includegraphics{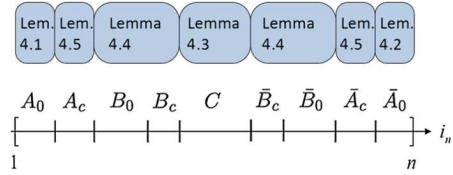}

\caption{Diagram of the sets of ranks as defined in (\protect\ref
{eqsetsranks}) and the
corresponding lemmas in Section~\protect\ref{secexpressions} which provide
expressions for
the local levels $\alpha_{i_n,n}$ for the different growth rates of $i_n$.}
\label{figseverallemmas}
\end{figure}

By means of Theorem~\ref{lemloclevPois} and the Stirling formula
\begin{equation}
\label{stirling}
i!=\sqrt{2 \uppi} i^{i+1/2} \exp(-i) \bigl[ 1+\mathrm{O} (1/i ) \bigr]
\end{equation}
as well as Theorem~\ref{zoneofconv}, we are now able to calculate
local levels $\alpha_{i,n}$ for various ranks $i$.
Local levels $\alpha_{i,n}$ of HC tests with critical values
$h_{i,n}(d_n(t))$ are given
in Lemmas \ref{loclevelsextremeintermediateleft}--\ref{loclevelnurun}.
For the sake of simplicity, we introduce the following notation for the
different
growth rates of $i_n$. We define the following sets of ranks $i_n\leq
n$, $n\in\mathbb{N}$,
\begin{eqnarray}
A_c &\widehat{=}&
i_n/u_n \to c \qquad\mbox{as } n\to\infty,
\nonumber\\
B_0 &\widehat{=}& i_n/u_n \to\infty\quad\mbox{and}\quad i_n/u_n^3 \to0\qquad \mbox{as } n\to\infty,
\nonumber\\
B_c &\widehat{=}& i_n/u_n^3 \to
c > 0 \qquad\mbox{as } n\to\infty,
\nonumber\\
\label{eqsetsranks}
C &\widehat{=}& i_n (1-i_n/n) /u_n^3
\to\infty \qquad\mbox{as } n\to \infty,
\\
\bar{B}_c &\widehat{=}& (n-i_n)/u_n^3
\to c > 0\qquad \mbox{as } n\to \infty,
\nonumber\\
\bar{B}_0 &\widehat{=}& (n-i_n)/u_n \to
\infty\quad\mbox{and}\quad (n-i_n)/u_n^3 \to0\qquad \mbox{as } n\to\infty,
\nonumber\\
\bar{A}_c &\widehat{=}& (n-i_n)/u_n \to c
\qquad\mbox{as } n\to\infty.
\nonumber
\end{eqnarray}
For example, for a sequence of ranks $i_n$, $n\in\mathbb{N}$, corresponding
to $A_c$ with $c=0$ we write $i_n\in A_0$.
Figure~\ref{figseverallemmas} summarizes which ranks $i_n$ correspond
to each lemma.
In the next two lemmas, we state local levels of HC tests for extreme
ranks and intermediate
ranks close to extreme ones, that is, $i_n \in A_0\cup \bar{A}_0$.

\begin{lemma}\label{loclevelsextremeintermediateleft}
For $i_n\in A_0$, we obtain
\begin{equation}
\label{alphasmallerintermediateleft}
\alpha_{i_n,n} = \frac{1}{\sqrt{2 \uppi i_n}} \biggl( \gamma
\frac{i_n}{u_n} \biggr)^{i_n} \exp ( - i_n v_n
) \biggl[ 1 + \mathrm{O} \biggl( \frac{1}{i_n} \biggr) +\mathrm{o}(1) \biggr],
\end{equation}
where $\gamma=\exp(1)/2$ and
\begin{equation}
\label{vn}
v_n = \bigl(\log(u_n)+2t-\log(
\uppi)+3i_n\bigr)/( 2 u_n) \bigl[ 1 + \mathrm{o}(1) \bigr].
\end{equation}
Alternatively, for $i_n \in A_0$ such that $i_n=\mathrm{o}(\sqrt{u_n})$ we get
\begin{equation}
\label{alphaextremeleft}
\alpha_{i_n,n} = \biggl( \frac{i_n^2}{2 u_n}
\biggr)^{i_n} \frac{1}{i_n!} \bigl[ 1 + \mathrm{o}( 1 ) \bigr]. 
\end{equation}
\end{lemma}

\begin{lemma}\label{loclevelsextremeintermediateright}
For $i_n\in\bar{A}_0$ we obtain
\begin{eqnarray}
\alpha_{i_n,n} & = &
\frac{\sqrt{\uppi}\exp(-2 t)}{\sqrt{2 (n-i_n)}}
\frac{1}{u_n \log(n)^2} \biggl( \frac{u_n}{ \gamma(n-i_n)} \biggr)^{n-i_n} \exp \bigl(
(n-i_n) w_n \bigr)
\nonumber
\\[-8pt]
\label{alphasmallerintermediateright}
\\[-8pt]
\nonumber
&&{}\times \bigl[ 1 + \mathrm{O} \bigl( 1/(n-i_n) \bigr) +\mathrm{o}(1) \bigr],
\end{eqnarray}
where $\gamma=\exp(1)/2$ and
\begin{equation}
\label{wn}
w_n = \bigl(\log(u_n)+2t-\log(
\uppi)+3(n-i_n)\bigr)/(2 u_n) \bigl[ 1 + \mathrm{o}(1) \bigr].
\end{equation}
Moreover, for $i_n \in\bar{A}_0$ such that $n-i_n=\mathrm{o}(\sqrt{u_n})$ we get
\begin{equation}
\label{alphaextremeright}
\alpha_{i_n,n} = \frac{\uppi\exp(-2t)}{u_n \log(n)^2} \biggl(
\frac{2u_n}{\exp(2)} \biggr)^{n-i_n} \frac{1}{(n-i_n)!} \bigl[ 1 + \mathrm{o}( 1 )
\bigr].
\end{equation}
\end{lemma}

The following lemma contains an expression for local levels of HC tests
for central ranks
and intermediates close to central ranks, that is, $i_n\in C$.

\begin{lemma}\label{loclevelslargerintermediatecentral}
Let $i_n\in C$. Then
\begin{equation}
\label{alphalargerintermediatecentral}
\alpha_{i_n,n} = \frac{\exp(-t)}{2u_n \log(n)} \biggl[ 1+ \mathrm{O} \biggl(
\frac
{\log(u_n)}{u_n} + \frac{u_n^{3/2}}{\sqrt{i_n(1-i_n/n)}} \biggr) \biggr],
\end{equation}
that is, local levels $\alpha_{i_n,n}$ with aforementioned $i_n$-values are
asymptotically equal.
Moreover, for a sequence $k_n\in\{ 1,\ldots,n \}$ such that
$k_n(1-k_n/n)/u_n^3\to\infty$
as $n\to\infty$ and all $i_n=k_n,\ldots,n-k_n$, local levels
$\alpha_{i_n,n}$ converge uniformly.
\end{lemma}

The next lemma provides local levels of HC tests corresponding to
intermediate ranks
$i_n \in  B_0 \cup B_c, \bar{B}_0 \cup\bar{B}_c$.

\begin{lemma}\label{un3}
Let $i_n \in B_0 \cup B_c$ or $i_n\in\bar{B}_0 \cup\bar{B}_c$.
Then
\begin{equation}
\label{alphaunpower3}
\alpha_{i_n,n} = \frac{\exp(-t) }{2 u_n \log(n)} \exp \biggl(
\frac{\sqrt{2} \zeta_n}{3} \bigl(u_n+ \mathrm{o}(u_n)\bigr) \biggr),
\end{equation}
where $\zeta_n = \sqrt{ u_n / i_n }$ if $i_n \in B_0 \cup B_c$ and
$\zeta_n = - \sqrt{ u_n / (n-i_n) }$ if $i_n \in\bar{B}_0 \cup\bar{B}_c$.
\end{lemma}

Finally, we give representations for local levels of HC tests for the remaining
intermediate ranks $i_n \in A_c$ and $i_n \in\bar{A}_c$.

\begin{lemma}\label{loclevelnurun}
Let $i_n \in A_c$ for a $c>0$ and set $c_n\equiv u_n/i_n$. Then
\begin{eqnarray}
\alpha_{i_n,n} & = & \frac{\sqrt{c_n}}{(1-\delta(c_n))\sqrt{2\uppi u_n}} \biggl[
\frac{\delta(c_n)\exp(1)}{\exp(\delta(c_n))} \biggr]^{{u_n}/{c_n}}
\nonumber
\\[-8pt]
\label{alphanurunleft}
\\[-8pt]
\nonumber
&&{}\times \bigl[ \sqrt{\uppi}\exp(-t) / \sqrt{u_n}
\bigr]^{({1-\delta(c_n)})/{\sqrt{c_n^2+2c_n}}} \bigl[1+\mathrm{o}(1)\bigr]
\end{eqnarray}
with $\delta(c)=1+c-\sqrt{c^2+2 c}$.
If $i_n \in\bar{A}_c$, $c>0$, and $c_n\equiv(n-i_n)/u_n$, then
\begin{eqnarray}
\alpha_{i_n,n} & = &
\biggl( 1 +
\frac{c_n}{1+\sqrt{1+2 c_n}} \biggr) \frac{ 1 }{\sqrt{2
\uppi c_n u_n}} \biggl( \frac{1+c_n+\sqrt{1+2 c_n}}{c_n }
\biggr)^{c_n u_n}
\nonumber\\
\label{alphanurunright}
&&{}\times\bigl( \sqrt{\uppi}\exp(-t)/ \sqrt{u_n}
\bigr)^{1+{1}/{\sqrt{1+2 c_n}}} \log(n)^{-(1+\sqrt{1+2 c_n})}
\\
&&
{}\times\bigl[1+\mathrm{o}(1)\bigr].\nonumber
\end{eqnarray}
\end{lemma}

\section{Monotonicity of HC local levels and related results} \label{secquotient}
First, we briefly illustrate the behavior of one-sided local levels of
HC tests for finite $n$-values.
\begin{figure}

\includegraphics{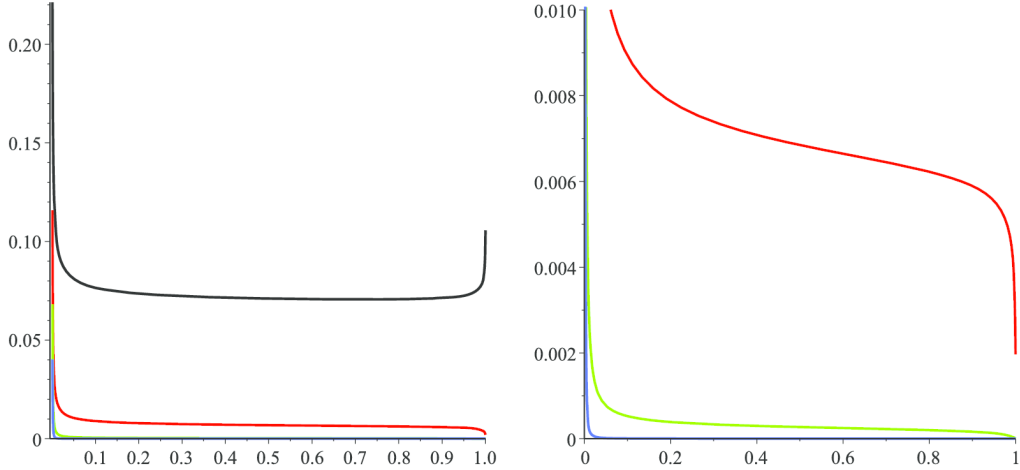}

\caption{The left graph: exact local levels curves calculated for
one-sided HC-tests
with $n=1000$ and $d=1.5,2.5,3.5,4.736$ (from top to bottom). The
right graph is zoomed.}
\label{figlocallevels}
\end{figure}
Figure~\ref{figlocallevels} provides exactly calculated local levels
$\alpha_{i,n}=\mathbb{P}(U_{i:n} < h_{i,n}(d))$ of HC tests $\varphi
^{\mathrm{HC}}$
(say) with critical values
$h_{i,n}(d)$, $i=1,\ldots,n$, for $n=1000$ and $d=1.5,2.5,3.5,4.736$.
For $d=1.5,2.5,3.5,4.736$ we get $\mathbb{E}_0(\varphi^{\mathrm{HC}}) =
0.803,0.322,0.111,0.05$, respectively.
That is, the HC test based on $d=4.736$ is a level $\alpha$ GOF test
for $\alpha=0.05$.
Figure~\ref{figlocallevels} illustrates that local levels are
decreasing for larger $d$-values.
Noting that our asymptotic investigations are given for $d\equiv d_n$
tending to infinity,
it seems that asymptotic results related to HC tests should be in
accordance with the
corresponding finite results for larger values of $d$.

Indeed, the following theorem shows that local levels $\alpha_{i,n}$ of
a HC test
with critical values $h_{i,n}(d_n(t))$, $i=1,\ldots,n$, are
asymptotically ($n\to\infty$)
non-increasing in $i$ in the following sense.
For non-decreasing sequences $i_n^{(1)}$ and
$i_n^{(2)}$ fulfilling $i_n^{(1)} < i_n^{(2)}$ for all $n\in\mathbb
{N}$, we get
$\lim_{n\to\infty}\alpha_{i_n^{(2)},n}/\alpha_{i_n^{(1)},n}\leq1$.
More precisely, $\alpha_{i_n^{(2)},n}/\alpha_{i_n^{(1)},n}$ depends on
the difference
$i_n^{(2)}-i_n^{(1)}$ and/or the ratio $i_n^{(1)}/i_n^{(2)}$.
Typically, the larger the difference $i_n^{(2)}-i_n^{(1)}$, the smaller
the ratio
$\alpha_{i_n^{(2)},n}/\alpha_{i_n^{(1)},n}$.

\begin{theorem}[(Asymptotic monotonicity of HC local levels)]\label{locallevelscomparisons1}
Let $i_n^{(1)} $ and $i_n^{(2)}$ be non-decreasing sequences that satisfy
$i_n^{(1)} < i_n^{(2)}$ for all $n\in\mathbb{N}$.
Let $\alpha_{i,n}$ denote the $i$th local level corresponding to a HC test
with critical values $h_{i,n}(d_n(t))$, $i=1,\ldots,n$.
Then
\begin{equation}
\label{aim1}
\lim_{n\to\infty}\alpha_{i_n^{(2)},n}/\alpha_{i_n^{(1)},n} = 1
\end{equation}
if the tuple $(i_n^{(1)}, i_n^{(2)})$ satisfies:
\begin{longlist}[(iii)]
\item[(i)] \vspace*{3pt}$(i_n^{(1)}, i_n^{(2)}) \in C\times C$,
\item[(ii)]  \vspace*{3pt}$(i_n^{(1)}, i_n^{(2)}) \in B_c \times B_c$,
\item[(iii)]  \vspace*{3pt}$(i_n^{(1)}, i_n^{(2)}) \in\bar{B}_c \times\bar{B}_c$,
\item[(iv)] $(i_n^{(1)}, i_n^{(2)}) \in B_0 \times B_0$
and  $i_n^{(1)}/i_n^{(2)} =1+ \mathrm{o} ( \sqrt{ i_n^{(1)}
/u_n^3}
)$,
\item[(v)]  $(i_n^{(1)}, i_n^{(2)}) \in\bar{B}_0 \times\bar{B}_0$ and $(n-i_n^{(2)})/(n-i_n^{(1)})
= 1+ \mathrm{o} ( \sqrt{ (n-i_n^{(2)}) /u_n^3} )$.
\end{longlist}
Moreover, we have
\begin{equation}
\label{aim01} 0 < \lim_{n\to\infty}\alpha_{i_n^{(2)},n} /
\alpha_{i_n^{(1)},n} < 1
\end{equation}
if one of the following conditions is satisfied:
\begin{longlist}[(viii)]
\item[(vi)]   \vspace*{3pt}$(i_n^{(1)}, i_n^{(2)}) \in B_{c_1} \times( B_{c_2} \cup C
\cup\bar{B}_c)$  with  $c_1 < c_2$,
\item[(vii)]   \vspace*{3pt}$(i_n^{(1)}, i_n^{(2)}) \in( C \cup\bar{B}_{c_1} )
\times\bar{B}_{c_2}$ with  $c_2 < c_1$,
\item[(vi)]   \vspace*{3pt}$(i_n^{(1)}, i_n^{(2)}) \in B_{c_1} \times( B_{c_2} \cup C
\cup\bar{B}_c)$  with  $c_1 < c_2$,
\item[(vii)]  \vspace*{1pt}$(i_n^{(1)}, i_n^{(2)}) \in( C \cup\bar{B}_{c_1} )
\times\bar{B}_{c_2}$  with  $c_2 < c_1$,
\item[(viii)]  \vspace*{1pt}$(i_n^{(1)}, i_n^{(2)}) \in B_0 \times B_0$  and
$i_n^{(1)} / i_n^{(2)} = 1 - c_n \sqrt{ i_n^{(1)} /u_n^3}$  with
$c_n \to c >0$,
\item[(ix)]  \vspace*{3pt}$(i_n^{(1)}, i_n^{(2)}) \in\bar{B}_0 \times\bar{B}_0$
 and  $(n-i_n^{(2)})/(n-i_n^{(1)}) = 1 - c_n \sqrt{ n-i_n^{(2)}
/u_n^3}$, with  $c_n \to c >0$,
\item[(x)]  \vspace*{3pt}$(i_n^{(1)}, i_n^{(2)}) \in A_c \times A_c$ with  $c>0$
 and  $i_n^{(2)} - i_n^{(1)} \equiv m$   for an  $m\in\mathbb{N}$,
\item[(xi)]  \vspace*{3pt}$(i_n^{(1)}, i_n^{(2)}) \in\bar{A}_c \times\bar{A}_c$
 with  $c>0$  and  $i_n^{(2)} - i_n^{(1)} \equiv m$  for an  $m\in
\mathbb{N}$.
\end{longlist}
Finally,
\begin{equation}
\label{aim0}
\lim_{n\to\infty}\alpha_{i_n^{(2)},n} /
\alpha_{i_n^{(1)},n} = 0
\end{equation}
for all other tuples with $ i_n^{(1)} < i_n^{(2)} $ when this limit exists.
\end{theorem}

\begin{figure}

\includegraphics{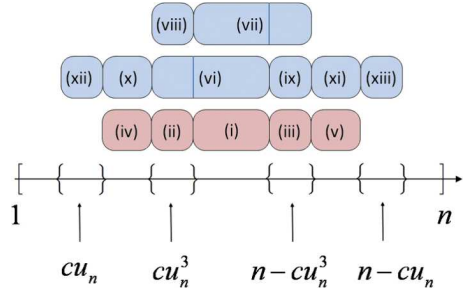}

\caption{Diagram of the sets of ranks as defined in (\protect\ref{eqsetsranks}) and the
corresponding regions covered by Theorem~\protect\ref
{locallevelscomparisons1}.}
\label{figTh4.1}
\end{figure}

Figure~\ref{figTh4.1} illustrates the regions of validity of (i)$-$(xi)
in Theorem~\ref{locallevelscomparisons1}.
Since $\alpha_{1,n}\to0$ as $n\to\infty$, cf. (\ref
{alphaextremeleft}), and local
levels $\alpha_{i,n}$, $i=2,\ldots,n$ are smaller than $\alpha_{1,n}$
for $n$ large enough,
cf. Theorem~\ref{locallevelscomparisons1}, the following result is obvious.

\begin{theorem}\label{allloclevelstendstozero}
For the local levels of the HC test it holds that
\[
\lim_{n\to\infty} \max_{1\leq i\leq n} \alpha_{i,n}
= 0.
\]
\end{theorem}

Theorem~\ref{allloclevelstendstozero} implies that local levels
corresponding to a
HC test show a completely different limiting behavior than the local
levels corresponding
to KS tests, cf. Figure~\ref{figlocallevelsKS}. Moreover, the
statement of
Theorem~\ref{allloclevelstendstozero} on the local levels of HC
tests vanishing asymptotically
allows us to deduce a result on the more general case of asymptotic
level $\alpha$ GOF tests
with prespecified local levels $\alpha_{i,n}$, $i=1,\ldots,n$.

\begin{remark}
For a level $\alpha$ GOF test with $\alpha_{i,n}$ satisfying (\ref
{deflocallevels1}) or (\ref{deflocallevels2}), we get
\[
\lim_{n\to\infty}\min_{1\leq i \leq n} \alpha_{i,n} =
0 \quad \mbox{or}\quad \lim_{n\to\infty} \min_{1\leq i \leq n}
\alpha_{i,n}^= = 0,
\]
respectively.
Thus, it is impossible to construct an asymptotic level $\alpha$ GOF
test with local levels which are all asymptotically bounded away from zero.
\end{remark}

Lemmas \ref{loclevelsextremeintermediateleft}--\ref
{loclevelnurun} and Theorem~\ref{locallevelscomparisons1}
lead to the next lemma that provides the asymptotics of level $\alpha$
GOF tests with equal local levels.
%

\begin{lemma}\label{equalloclevasymptnew}
For one- or two-sided GOF tests with local levels equal to $\alpha
_n^{\mathrm{loc}}$, $n\in\mathbb{N}$, we obtain an asymptotic level
$\alpha$ test iff
\[
\lim_{n\to\infty}\alpha_n^{\mathrm{loc}} \cdot
\frac{2 \log(\log
(n)) \log(n) }{-\log
(1-\alpha)} = 1.
\]
\end{lemma}

The rather technical and straightforward proof will be presented in a
forthcoming paper.

\begin{remark}
Lemma~\ref{equalloclevasymptnew} is up to now the most precise result
concerning the asymptotics of the test with equal local levels.
For example, adapting Theorem~4.1 in the third version of \cite
{moscovich} leads to an asymptotic interval for $\alpha_n^{\mathrm{loc}}$.
Moreover, results in \cite{kaplan} and \cite{kaplannew} can be seen as
a very rough approximation for the rate given in Lemma~\ref
{equalloclevasymptnew}.
\end{remark}

\section{Comparison of GOF tests in the finite case} \label{secGOFcomparison}
In this section, we compare one-sided versions of KS tests $\varphi
^{\mathrm{KS}}$, HC tests $\varphi^{\mathrm{HC}}$ and GOF tests $\varphi(\alpha
_n^{\mathrm{loc}})$ with equal local levels for a finite sample size $n$.
In order to compare these tests in a fair way, all considered tests
will be of exact level $\alpha$.
That is, for fixed $n\in\mathbb{N}$ and $\alpha\in(0,1)$ we determine
parameters of the considered tests, that is, find $c$ for the KS test with
critical values $i/n-c/\sqrt{n}$, $i=1,\ldots,n$, a parameter $d$ for
the HC test based on $h_{i,n}(d)$, $i=1,\ldots,n$, given in~(\ref{hn})
and $\alpha_n^{\mathrm{loc}}$ for the GOF test with equal local
levels, so that
\[
\mathbb{E}_0\bigl(\varphi^{\mathrm{KS}}\bigr) =
\mathbb{E}_0\bigl(\varphi^{\mathrm{HC}}\bigr) =
\mathbb{E}_0\bigl(\varphi\bigl(\alpha _n^{\mathrm{loc}}
\bigr)\bigr)=\alpha.
\]
Clearly, such parameters can be found numerically,
for example, via some search algorithm, whenever the probability to
reject the true null hypothesis can be numerically calculated.
Thereby, the computation of the joint c.d.f. of the order statistics
$U_{1:n},\ldots,U_{n:n}$, that is, $\mathbb{P}(U_{i:n} \leq c_i,
i=1,\ldots,n)$, plays
the key role in the one-sided case, while the computation of $\mathbb
{P}(c_i <
U_{i:n} < \tilde{c}_i, i=1,\ldots,n)$ is crucial in the two-sided case.
Probabilities of the first type can be calculated by Noe's, Bolshev's,
Steck's or Khmaladze's recursions, $\mathbb{P}(c_i < U_{i:n} < \tilde{c}_i,
i=1,\ldots,n)$ can be calculated by Noe's, Ruben's or Khmaladze's recursions,
for example, cf. \cite{khmaladze} and  pages 357--370 in \cite{shorackwellner}.
If the sample size $n$ is so large that exact computations are no
longer possible, that is, $n\gg10^4$, the parameters $d$ and $\alpha
_n^{\mathrm{loc}}$ can approximately be calculated via numerical simulations.

For example, for $\alpha=0.05$ and $n=100$, $500$, $1000$ we get\vspace*{1pt} by
numerical calculations $\mathbb{E}_0(\varphi^{\mathrm{HC}}) = \alpha$ for
$d=4.725$,
$4.734$, $4.736$, respectively, $\mathbb{E}_0(\varphi^{\mathrm{KS}}) =
\alpha$ for
$c=1.22387$ and $\mathbb{E}_0(\varphi(\alpha_n^{\mathrm
{loc}}))=\alpha$ for $\alpha
_n^{\mathrm{loc}}=0.00246,0.00145,0.00122$, respectively.
The asymptotic local level in Lemma~\ref{equalloclevasymptnew} is equal
to $0.00365,0.00226$, $0.00192$ for $\alpha=0.05$ and $n=100,500,1000$,
respectively, so that the asymptotic local level seems to be larger
than the finite counterpart $\alpha_n^{\mathrm{loc}}$.

\begin{figure}

\includegraphics{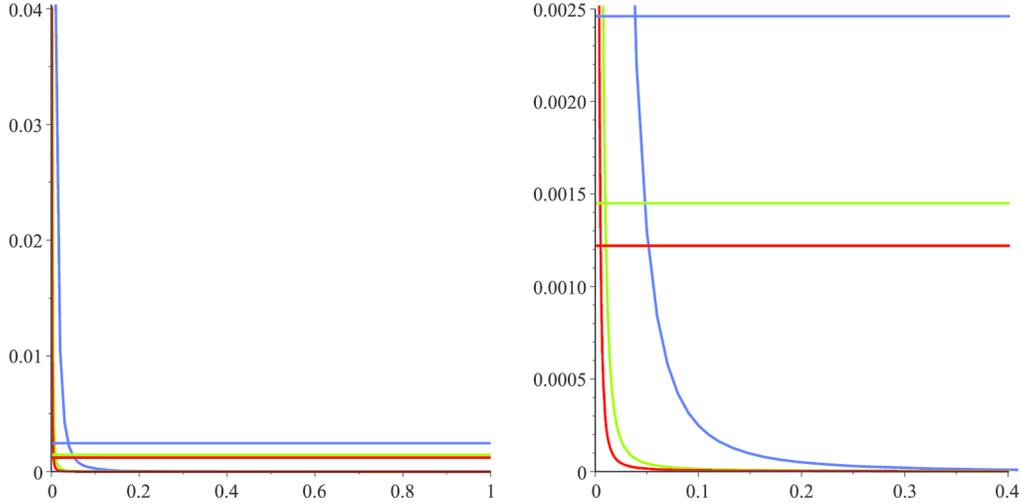}

\caption{Local levels curves corresponding to HC tests $\varphi
^{\mathrm{HC}}$ based on $h_{i,n}(d)$ with $d=4.725, 4.734,4.736$
(curves from top
to bottom) for $n=100,500,1000$, respectively, leading to $\mathbb
{E}_0(\varphi
^{\mathrm{HC}}) = 0.05$, and the corresponding local levels $\alpha
_n^{\mathrm{loc}}=0.00246,0.00145,0.00122$
(straight lines from top to bottom) that imply $\mathbb{E}_0(\varphi
(\alpha
_n^{\mathrm{loc}})) = 0.05$ for GOF tests with local levels equal
$\alpha
_n^{\mathrm{loc}}$.
The right graph is zoomed.}
\label{figloclevelsHCandGOFn1005001000}
\end{figure}

Figure~\ref{figloclevelsHCandGOFn1005001000} shows local
levels curves of the level $\alpha$ HC tests together with equal local
levels $\alpha_n^{\mathrm{loc}}$ (straight lines) for $n=100,500,1000$.
Local levels of the corresponding KS tests are given in Figure~\ref{figlocallevelsKS}.
Note that almost all local levels of the HC tests are smaller than the
corresponding $\alpha_n^{\mathrm{loc}}$ and only the first ones are larger,
for example, for $n=100,500,1000$ we get $\alpha_{i,n} \geq\alpha
_n^{\mathrm{loc}}$ if $i\leq3,4,5$, respectively, and $\alpha_{i,n}
< \alpha
_n^{\mathrm{loc}}$ else.
This indicates higher sensitivity of the GOF test with equal local
levels in a specific intermediate range than by the HC tests.

Now we consider the aforementioned level $\alpha$ GOF tests in terms of
their rejection curves.
\begin{figure}

\includegraphics{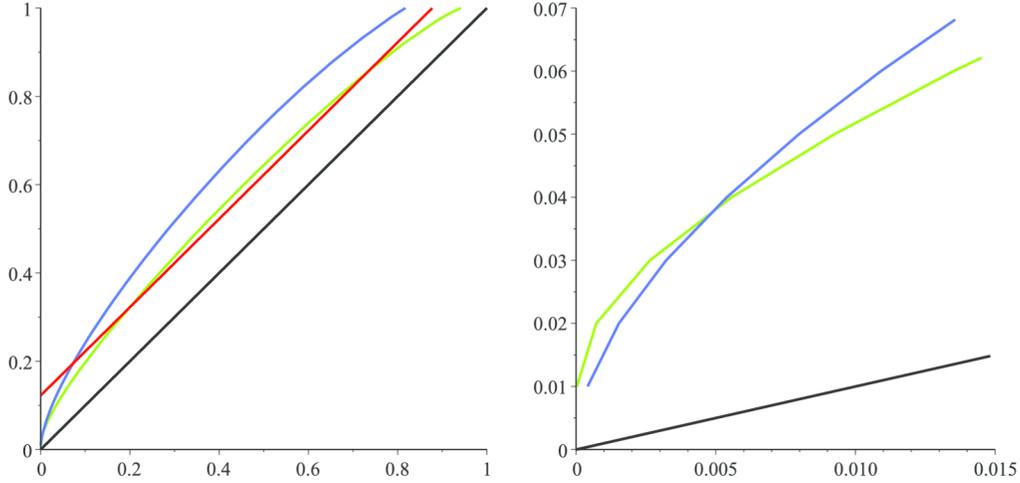}

\caption{Rejection curves of the level $\alpha$ GOF tests
$\varphi^{\mathrm{HC}}$,
$\varphi(\alpha_n^{\mathrm{loc}})$ and $\varphi^{\mathrm{KS}}$
together with the
diagonal (from top to bottom in $0.5$, respectively) for $n = 100$,
$\alpha=0.05$, $\varphi^{\mathrm{HC}}$
based
on
$h_{i,n}(d)$, $i=1,\ldots,n$,
with $d = 4.725$, $\varphi(\alpha_n^{\mathrm{loc}})$ based on
$\alpha_n^{\mathrm{loc}}=0.00246$ and $\varphi^{\mathrm{KS}}$
based on $i/n-c/\sqrt{n}$, $i=1,\ldots,n$
 with $c = 1.22387$. The right graph is zoomed.}
\label{figloclevelsallequal}
\end{figure}
Figure~\ref{figloclevelsallequal} shows rejection curves for $n=100$.
Here, critical values induced by $\varphi(\alpha^{\mathrm{loc}}_n)$
are larger than the corresponding HC critical values for $i\geq4$ and
only slightly smaller than the KS critical values in a specific central
range, while the latter are considerably smaller in tails.
Moreover, although all considered tests are level $\alpha$ tests,
almost all of the HC critical values are considerably smaller than the
corresponding critical values of the GOF test $\varphi(\alpha
_n^{\mathrm{loc}})$ with equal local levels.
It indicates that the smallest critical values have the biggest impact
on $\mathbb{E}_0(\varphi)$ for any GOF test $\varphi$ while other critical
values influence $\mathbb{E}_0(\varphi)$ only slightly.
Further exact calculation showed that a similar picture is observed for
various $n$-values.

Altogether, it seems that the level $\alpha$ GOF tests with equal local
levels offer a good alternative to the classical GOF tests especially
if it is not clear what kind of deviation from the null hypothesis may occur.
For power comparisons between GOF tests with equal local levels and
other GOF tests see \cite{aldor13,GLF2014} and \cite{kaplannew}.

\section{Concluding remarks}\label{outlook}
In this paper, we introduced the concept of local levels $\alpha_{i,n}$
for a certain class of GOF tests.
These quantities serve as an indicator of regions of high/low local
sensitivity of a test and thus provide a method to compare tests with
respect to areas of sensitivity.
For example, the classical KS test has higher power for alternatives
that differ from the null
distribution in the central range.
This coincides with the fact that local levels of the KS tests are
considerably larger in the central range and are even equal to zero for
extremes and smaller intermediates.
In high-dimensional data with only sparse signals that are to be
detected, it would be advisable to perform a GOF test (or related
multiple tests) which is sensitive in the tails.
In such situations performing HC tests, which are asymptotically
sensitive only in the moderate tails, would be an advantage.
Due to the fact that the number of local levels corresponding to
central ranks is considerably higher than the number of local levels
corresponding to intermediate ranks, one may guess that the HC local
levels $\alpha_{i_n,n}$ for central ranks are much smaller than their
counterparts in the moderate tails.
Therefore, it is a rather striking result that central local levels are
indeed asymptotically as large as the ones in the moderate tails.
The reason for this may be hidden in the complex dependence structure
of order statistics, so that a further investigation in this direction
is needed.
In general, it seems to be an interesting issue to analyze local levels
of other multiple testing related GOF tests, thus gaining deeper
insight into their nature.
Figure~\ref{figloclevphidiv} suggests that the Berk--Jones test
comes close to the equal local levels test.
It might be of interest to compare the asymptotic
local levels of these tests as outlined in this paper for the HC test.
An additional difficulty is that explicit critical values needed in
(\ref{defGOF}) and (\ref{defGOF2}) are hard
to obtain for the Berk--Jones
test as well as for most of the other phi-divergence tests.

Furthermore, the concept of local levels may be used to construct new
tailored GOF tests if one has an idea in which region, that is, for
which kinds of alternatives, a test needs to be sensitive.
Given a set of suitable local levels we illustrated a way how to
construct the corresponding GOF test.
Moreover, by means of results related to the HC tests we showed that
there is no level $\alpha$ GOF tests with local levels asymptotically
uniformly bounded away from zero.
In view of the fact that most of the HC local levels are asymptotically
equal and that the first HC local levels are much too large so that the
remaining ones are too small in the finite case, the GOF test with
equal local levels $\alpha_{i,n} \equiv\alpha_n^{\mathrm{loc}}$
seems to be
a good alternative for the classical HC test, which is known for its
extremely slow asymptotics.
Although we do not have any explicit formula for the local level
$\alpha
_n^{\mathrm{loc}}$ as a function of the sample size $n$ and predefined level
$\alpha$,
we provide an asymptotic rate for $\alpha_n^{\mathrm{loc}}$ leading
to the
asymptotic level $\alpha$ test.\vspace*{-3pt}

\begin{appendix}
\section{Proofs of Sections \texorpdfstring{\protect\ref{secHC}}{2} 
and \texorpdfstring{\protect\ref{secapproximations}}{3}}\label{AppA}
\begin{pf*}{Proof of Lemma~\protect\ref{entwicklungh}}
Setting $A_n = 4 (d_n^2+n) i_n^2 /(n(d_n^2+2i_n)^2) $, a critical value
$h_{i_n,n}\equiv h_{i_n,n}(d_n(t)) $ can be represented as
\begin{equation}
\label{hformel}
h_{i_n,n} = \bigl(d_n^2 + 2
i_n\bigr) (1-\sqrt{1- A_n}) /\bigl(2 \bigl(
d_n^2+n\bigr) \bigr).\vspace*{-15pt}
\end{equation}
\begin{longlist}
\item[(i)] Let $i_n$ be such that $i_n=\mathrm{o}(u_n)$.
Since $A_n = \mathrm{O}(i_n^2/u_n^2)$, $A_n\to0$ as $n\to\infty$.
Applying the Taylor series $ 1-\sqrt{1-x} = x/2 + \mathrm{O} (x^2) $ for $x\in
(0,1)$, we get
\begin{eqnarray*}
h_{i_n,n} & = & \bigl(d_n^2 + 2 i_n
\bigr)/ \bigl(2 \bigl( d_n^2+n\bigr)\bigr) \bigl[ 2
\bigl(d_n^2+n\bigr) i_n^2 / n /
\bigl(d_n^2+2i_n\bigr)^2 + \mathrm{O}
\bigl( i_n^4 / u_n^4 \bigr) \bigr]
\\
& = & i_n^2 /\bigl( n
\bigl(d_n^2+2i_n\bigr)\bigr) + \mathrm{O} \bigl(
i_n^4 / \bigl(n u_n^3\bigr)
\bigr)
\nonumber
\end{eqnarray*}
and hence $ n u_n h_{i_n,n} / i_n^2 = u_n/( d_n^2+2i_n) + \mathrm{O}  (
i_n^2/u_n^2  ) $.
Noting that
\begin{equation}
\label{dformel} d_n^2 = 2 u_n +
\log(u_n) + 2t -\log(\uppi) + \mathrm{O} \bigl( \log(u_n)^2
/ u_n \bigr)
\end{equation}
and $1/(2+x)=1/2-x/4+\mathrm{O}(x^2)$ for $x\in(0,1)$, we get
\[
\frac{u_n}{ d_n^2+2i_n} = \frac{1}{2} - \frac{\log(u_n)+2 t-\log(\uppi)+2 i_n}{4 u_n} + \mathrm{O} \biggl(
\frac
{\log(u_n)^2+ i_n^2}{u_n^2} \biggr)
\]
and consequently (\ref{hseriesextremeleft}) follows.

\item[(ii)]  Let $c_n\equiv u_n/i_n\to c$ as $n\to\infty$ for some fixed $c>0$.
Obviously, $A_n = [(d_n^2+2i_n)/(2i_n)]^{-2}(1+\mathrm{O}(u_n/n))$.
Since $(1+x)^{-2}=1-2x+\mathrm{O}(x^2)$ for $x\in(0,1)$, (\ref{dformel}) leads to
\begin{eqnarray*}
\biggl[ \frac{d_n^2+2i_n}{2i_n} \biggr]^{-2} & = & \biggl[ (1+
c_n ) + \frac{\log(u_n)+2 t - \log(\uppi)}{2 i_n} +\mathrm{O} \biggl( \frac
{\log(u_n)^2}{u_n^2} \biggr)
\biggr]^{-2}
\\
& = & \frac{1}{(1+ c_n )^2} \biggl[ 1 - \frac{\log(u_n)+2 t - \log(\uppi)}{ i_n  (1+ c_n  ) } +\mathrm{O} \biggl(
\frac{\log(u_n)^2}{u_n^2} \biggr) \biggr].
\end{eqnarray*}
Since $ 1-\sqrt{1-a(1-x)}=1-\sqrt{1-a}- a x / (2\sqrt{1-a}) +\mathrm{O}(x^2) $
for $a>0$ and $x\in(0,1)$, we get
\[
1-\sqrt{1-A_n} = 1- \frac{\sqrt{ c_n^2+2c_n }}{1+c_n} - \frac{\log
(u_n)+2 t-\log(\uppi) }{2 i_n (1+c_n)^2 \sqrt{c_n^2+2c_n}} + \mathrm{O} \biggl( \frac{\log(u_n)^2}{u_n^2} \biggr).
\]
Furthermore,
\begin{eqnarray*}
\frac{d_n^2+2i_n}{2 (d_n^2+n)} & = & \bigl( d_n^2+2i_n
\bigr)/(2n) \bigl[ 1 + \mathrm{O} \bigl( d_n^2/ n \bigr) \bigr]
 =  (i_n / n) \bigl[ 1 + d_n^2 /(2
i_n)+ \mathrm{O} ( u_n/n ) \bigr]
\\
& = & (i_n/ n) \biggl[ 1 + c_n + \frac{ \log(u_n)+2 t-\log(\uppi)}{2 i_n} + \mathrm{O} \bigl( \log(u_n)^2 /u_n^2 \bigr)
\biggr].
\end{eqnarray*}
Formula (\ref{hformel}) immediately leads to
\begin{eqnarray*}
h_{i_n,n} & = & i_n \Bigl(1+c_n-\sqrt
{c_n^2 +2c_n}\Bigr)\big/n
\\
& & {}\times\Bigl[ 1 + \mathrm{O} \bigl( \log(u_n)^2/u_n^2 \bigr)- \bigl( \log(u_n)+2 t-\log(\uppi) \bigr)\big/ \Bigl(2 i_n
\sqrt{c_n^2+2c_n} \Bigr) \Bigr]
\end{eqnarray*}
and hence, we get (\ref{hseriesintermediate}).

\item[(iii), (iv)] Now, let $u_n=\mathrm{o}(i_n)$.
Due to $1/(1+x)^2=1-2x+3x^2+\mathrm{O}(x^3)$ for $x\in(0,1)$, we get
\begin{eqnarray*}
A_n & = & \bigl(d_n^2+n\bigr)/n \bigl[
\bigl(d_n^2+2i_n\bigr)/(2 i_n)
\bigr]^{-2}
 =  \biggl( 1 + \frac{d_n^2}{n} \biggr) \biggl[ 1 - \frac{d_n^2}{ i_n} +
\frac{3
d_n^4}{ 4 i_n^2} + \mathrm{O} \biggl( \frac{d_n^6}{ i_n^3} \biggr) \biggr]
\\
& = & 1 - \frac{d_n^2}{i_n} \biggl(1-\frac{i_n}{n} \biggr) +
\frac{d_n^4}{ i_n^2} \biggl(\frac{3}{4}-\frac{i_n}{n} \biggr) + \mathrm{O} \biggl( \frac{d_n^6}{ i_n^3} \biggr).
\end{eqnarray*}
Hence, for $i_n$ such that $u_n=\mathrm{o}(i_n(1-i_n/n))$ we arrive at
\begin{eqnarray*}
A_n & = & 1 - 2 u_n ( 1- i_n/ n )/
i_n - \bigl(\log(u_n)+2 t-\log(\uppi) \bigr) (
1-i_n/ n ) / i_n
\\
& &{}+ \mathrm{O} \bigl( u_n^2/ i_n^2 +
\log(u_n)^2 (1-i_n/n ) /(i_n
u_n) \bigr)
\end{eqnarray*}
and for $i_n$ such that $n-i_n=\mathrm{O}(u_n)$ we obtain
\begin{eqnarray*}
A_n & = & 1 - u_n^2 \bigl(1 +2 (
n-i_n)/ u_n \bigr) /i_n^2-
u_n \bigl(\log (u_n)+2t-\log(\uppi)\bigr) /
i_n^2
\\
& & {}\times \bigl(1 + (n-i_n)/u_n \bigr) + \mathrm{O} \bigl(
\log(u_n)^2 / i_n^2 \bigr).
\end{eqnarray*}
Then
\begin{eqnarray*}
1-\sqrt{1-A_n} & = & 1 - \sqrt{ \frac{2 u_n}{i_n} \biggl(1-
\frac{i_n}{n} \biggr) } - \frac
{\log
(u_n)+2 t-\log(\uppi)}{2\sqrt{2 i_n u_n}} \sqrt{ 1 - i_n/n }
\\
& & {}+ \mathrm{O} \biggl( \frac{u_n^{3/2}}{i_n^{3/2}\sqrt{1-i_n/n}} + \frac{ \log
(u_n)^2}{\sqrt{i_n} u_n^{3/2} } \sqrt{1-
i_n/ n} \biggr)
\end{eqnarray*}
for $i_n$ such that $u_n=\mathrm{o}(i_n(1-i_n/n))$ and
\begin{eqnarray*}
1-\sqrt{1-A_n} & = & 1 - \frac{u_n}{i_n}\sqrt{1+2
\frac{n-i_n}{u_n}} - \frac{\log(u_n)+2 t -\log(\uppi)}{2i_n}
\\
& & {}\times\frac{1+({n-i_n})/{u_n}}{\sqrt{1+2({n-i_n})/{u_n}}} + \mathrm{O} \biggl( \frac{\log(u_n)^2}{i_n u_n} \biggr)
\end{eqnarray*}
if $n-i_n=\mathrm{O}(u_n)$.
Since
\begin{eqnarray*}
\frac{d_n^2+2i_n}{2(d_n^2+n)} & = & \biggl( \frac{i_n}{n} + \frac{d_n^2}{2n} \biggr)
\biggl[ 1 - \frac{d_n^2}{n} + \mathrm{O} \biggl( \frac{u_n^2}{n^2} \biggr) \biggr]
\\
& = & \frac{i_n}{n} \biggl[ 1 + \biggl( 1 - \frac{2 i_n}{n} \biggr)
\frac
{2u_n+\log
(u_n)+2 t-\log(\uppi)}{2 i_n}
\\
&&\quad\hspace*{5pt}{}+ \mathrm{O} \bigl( u_n^2/(n i_n) + ( 1 - 2
i_n/ n ) \log(u_n)^2 / (u_n
i_n) \bigr) \biggr]
\end{eqnarray*}
for all $i_n\leq n$, we get for $i_n$ such that $u_n = \mathrm{o}(i_n(1-i_n/n))$
\begin{eqnarray}
 h_{i_n,n} & = &
\frac{i_n}{n} \biggl[ 1 -
\sqrt{\frac{2 u_n}{i_n} \biggl( 1- \frac{i_n}{n} \biggr) } -
\frac{\log(u_n)+2 t-\log(\uppi)}{2\sqrt{2 i_n u_n}}\sqrt{ 1- \frac
{i_n}{n} }
\nonumber
\\
\label{uniformlyconvergenceofh}
& & \hspace*{15pt}{}+ ( 1 - 2 i_n / n ) \bigl( 2u_n+
\log(u_n)+2 t-\log(\uppi) \bigr) / (2 i_n)
\\
&&
\nonumber
\hspace*{15pt}{}+ \mathrm{O} \biggl( \frac{u_n^{3/2}}{i_n^{3/2}\sqrt{1-i_n/n}} + \frac{ \log(u_n)^2}{\sqrt{i_n} u_n^{3/2} } \sqrt{1-
\frac
{i_n}{n}} \biggr) \biggr]
\end{eqnarray}
and
\begin{eqnarray*}
h_{i_n,n} & = & \frac{i_n}{n} \biggl[ 1 - \frac{u_n}{i_n} \bigl(
1 + \sqrt{1+2(n-i_n)/u_n } \bigr) -\frac{\log(u_n)+2t-\log(\uppi)}{2 i_n}
\\
& & \hspace*{15pt}{}\times \biggl( 1+ \frac{1+ (n-i_n)/u_n }{ \sqrt{1+2 (n-i_n)/u_n} } \biggr) + \mathrm{O} \biggl( \frac{\log(u_n)^2}{i_n u_n}
\biggr) \biggr]
\end{eqnarray*}
for $i_n$ with $n-i_n=\mathrm{O}(u_n)$.\hfill\qed
\end{longlist}
\noqed\end{pf*}

\begin{pf*}{Proof of Theorem~\protect\ref{lemloclevPois}}
We only prove (\ref{eqloclevPois})--(\ref{eqloclevPois2}).
The cases (\ref{eqloclevPoisright})--(\ref{eqloclevPois22}) can
be handled analogously.
For $n\geq4$, $h_{i_n,n}\leq1/4$ and $|i_n-n h_{i_n,n}|\leq\sqrt
{n}/2$ we obtain
\[
\frac{\mathbb{P}(X_n=i_n)}{\mathbb{P}(Y_n = i_n)} = 1 -\frac{(i_n-n h_{i_n,n})^2}{2 n (1-h_{i_n,n})} + \frac{i_n}{2 n} + \mathrm{O} \biggl(
\frac{(i_n-n h_{i_n,n})^3}{ n^2 } + \frac{i_n^2}{n^2} \biggr),
\]
cf.  formula (17) in \cite{Prokhorov}.
Therefore, (\ref{hseriesextremeleft})--(\ref{hseriescentralintermed})
in Lemma~\ref{entwicklungh} lead to
\begin{equation}
\label{loweralp}
\alpha_{i_n,n} = \mathbb{P}(Y_n \geq
i_n) \bigl[ 1 + \mathrm{o}(1) \bigr]\qquad \mbox{at least for } i_n=\mathrm{o} (
\sqrt{ n / u_n} ).
\end{equation}
Moreover, for $i_n=\mathrm{o} (\sqrt{n/u_n} )$ we get
\[
\alpha_{i_n,n} = \mathbb{P}(Y_n = i_n) \Biggl( 1+
\sum_{k=1}^\infty \prod
_{j=1}^{k} \frac{n h_{i_n,n}}{i_n+j} \Biggr) \bigl[1+\mathrm{o}(1)
\bigr],
\]
cf. \cite{shorackwellner}, page~485.
Since $n h_{i_n,n}/i_n <1$ for larger $n$-values, we get
\[
1 < 1+ \sum_{k=1}^\infty\prod
_{j=1}^{k} \frac{n h_{i_n,n}}{i_n+j} \leq\sum
_{k=0}^\infty \biggl( \frac{n h_{i_n,n}}{i_n+1}
\biggr)^k = \frac{i_n+1}{i_n-n h_{i_n,n}+1}
\]
for $i_n=\mathrm{o}(\sqrt{n/u_n})$.
Hence, $\lim_{n\to\infty}(i_n+1)/(i_n-n h_{i_n,n}+1) = 1$ in case
$i_n = \mathrm{o}(u_n)$,
that is, (\ref{eqloclevPois}) follows.
For $i_n$ such that $u_n/i_n\to c>0$, $n\to\infty$,
\[
\lim_{n\to\infty}\frac{i_n+1}{i_n-n h_{i_n,n}+1} = \frac{1}{\sqrt
{c^2+2 c}-c}
\]
and for a fixed $k\in\mathbb{N}$
\[
\lim_{n\to\infty} \Biggl( \prod_{j=1}^{k}
\frac{n h_{i_n,n}}{i_n+j} \Biggr) \Big/ \biggl( \frac{n h_{i_n,n}}{i_n+1} \biggr)^k =1,
\]
which implies (\ref{inCun}).
Furthermore, from known asymptotic decompositions for the incomplete
gamma function (e.g.,
cf. \cite{ivchenko,Ferreira} and \cite{Bateman}, page~140) and
from the fact that for
$i_n=\mathrm{o}(\sqrt{n/u_n})$ such that $u_n=\mathrm{o}(i_n)$ it holds $i_n\to\infty$,
$n h_{i_n,n}\to\infty$ and
$(i_n-n h_{i_n,n})/\sqrt{n h_{i_n,n}}\to\infty$, $n\to\infty$, we obtain
\[
\mathbb{P}(Y_n\geq i_n) = i_n/(i_n-n
h_{i_n,n}-1 ) \mathbb {P}(Y_n=i_n) \bigl[1+\mathrm{o}(1)
\bigr].
\]
This together with (\ref{hseriescentralintermed}) and (\ref
{loweralp}) imply (\ref{eqloclevPois2}).
\end{pf*}

\section{Proofs of Section~\texorpdfstring{\protect\ref{secexpressions}}{4}}\label{AppB}

\begin{pf*}{Proof of Lemma~\protect\ref{loclevelsextremeintermediateleft}}
With respect to Theorem~\ref{lemloclevPois}, it suffices to
calculate $\mathbb{P}(Y_n=i_n)$, where $Y_n\sim\mathcal{P}( n h_{i_n,n})$.
Obviously,
\[
\mathbb{P}(Y_n=i_n) = \bigl(i_n^{i_n}/
i_n!\bigr) \bigl[ ( n h_{i_n,n} / i_n ) \exp (- n
h_{i_n,n} / i_n ) \bigr]^{i_n}.
\]
Since (\ref{hseriesextremeleft}) implies $n h_{i_n,n}/i_n = \mathrm{o}(1)$
and $x \exp(-x)= x-x^2+\mathrm{O}(x^3)$ for $x\in(0,1)$, we obtain
\[
\mathbb{P}(Y_n=i_n) = \bigl(i_n^{i_n} /
i_n!\bigr) \bigl[ (n h_{i_n,n} / i_n) - ( n
h_{i_n,n} / i_n )^2 + \mathrm{O} ( n h_{i_n,n} /
i_n )^3 \bigr]^{i_n}.
\]
Setting representation (\ref{hseriesextremeleft}) for a critical
value $h_{i_n,n}$ in the equation above, we get
\begin{equation}
\label{fast}
\mathbb{P}(Y_n=i_n) = (B_n /
i_n!) \bigl( i_n^2 / 2 u_n
\bigr)^{i_n},
\end{equation}
where
\begin{eqnarray*}
B_n & = & \biggl( 1- \frac{\log(u_n)+2t-\log(\uppi)+3i_n}{2 u_n} + \mathrm{O} \biggl(
\frac{\log(u_n)^2+i_n^2}{ u_n^2} \biggr) \biggr)^{i_n}
\\
& = & \exp\bigl(i_n \log\bigl( 1- \bigl(\log(u_n)+2t-
\log(\uppi)+3i_n\bigr)/(2 u_n)
+ \mathrm{O}\bigl( \bigl(\log(u_n)^2+i_n^2
\bigr)/ u_n^2 \bigr)\bigr)\bigr).
\end{eqnarray*}
Since $i_n/u_n\to0$ as $n\to\infty$ and $\log(1-x)=-x+\mathrm{O}(x^2)$ as
$x\to
0$, it follows
\begin{equation}
\label{Bn} B_n = \exp \biggl( -i_n \biggl[
\frac{\log(u_n)+2t-\log(\uppi)+3 i_n}{2 u_n} + \mathrm{O} \biggl( \frac
{\log(u_n)^2+i_n^2}{ u_n^2} \biggr) \biggr] \biggr).
\end{equation}
Particularly, for $i_n=\mathrm{o}(\sqrt{u_n})$ we get $B_n = 1+\mathrm{o}(1) $, so that
(\ref{fast}) immediately leads to (\ref{alphaextremeleft}).
Finally, the Stirling formula (\ref{stirling}), (\ref{fast}) and
(\ref{Bn}) imply (\ref{alphasmallerintermediateleft}).
\end{pf*}

\begin{pf*}{Proof of Lemma~\protect\ref{loclevelsextremeintermediateright}}
Due to Theorem~\ref{lemloclevPois}, we have to calculate
\[
\mathbb{P}(\tilde{Y}_n = n-i_n) = (n-n
h_{i_n,n})^{n-i_n} \exp ( -n+n h_{i_n,n} ) /
(n-i_n)!.
\]
Setting $c_n\equiv(n-i_n)/u_n$, we get $c_n\to0$ as $n\to\infty$.
In order to simplify (\ref{hseriesextremeright}), we obtain $1+\sqrt
{1+2 c_n} = 2 + c_n - c_n^2/2+\mathrm{O}(c_n^3)$
and $1 + (1+c_n)/\sqrt{1+2c_n} = 2+\mathrm{O}(c_n^2)$. Then
\begin{eqnarray}
n h_{i_n,n} & = & i_n-u_n \bigl( 2
+c_n - c_n^2 / 2 \bigr) -
\log(u_n)-2t+\log(\uppi)
\nonumber\\
& &
{}+ \mathrm{O}\bigl( u_n c_n^3 +
\log(u_n)c_n^2+ \log(u_n)^2
/ u_n \bigr),\nonumber
\\
\exp ( -n+n h_{i_n,n} )
& = & \frac{\uppi\exp(-2t)}{u_n \log(n)^2 \exp(2(n-i_n))} \bigl[ 1 + \mathrm{O} \bigl( \log(u_n)^2/
u_n \bigr) \bigr]
\nonumber
\\[-8pt]
\label{expright}
\\[-8pt]
\nonumber
&&
{}\times\exp \bigl( (n-i_n) \bigl[ c_n/2 + \mathrm{O} \bigl( c_n^2 + c_n\log (u_n)/u_n
\bigr) \bigr] \bigr)
\end{eqnarray}
and
\begin{eqnarray*}
&&\!\!\! (n-n h_{i_n,n})^{n-i_n}\\
&&\!\!\!\hspace*{9pt} =
(2
u_n)^{n-i_n}\exp \bigl( (n-i_n)\log \bigl( (n-n
h_{i_n,n})/(2u_n) \bigr) \bigr)
\\
&&\!\!\!\hspace*{9pt}=
(2 u_n)^{n-i_n} \exp\biggl(
(n-i_n)\log\biggl(  1+c_n
+ \frac{\log(u_n)+2t-\log(\uppi)}{2 u_n} + \mathrm{O} \bigl( c_n^2 +
\log(u_n)^2/u_n^2 \bigr)
\biggr) \biggr).
\end{eqnarray*}
Taylor's series $\log(1+x)=x+\mathrm{O}(x^2)$ for $x\in( 0,1)$ leads to
\begin{eqnarray*}
&& (n-n h_{i_n,n})^{n-i_n}\\
&&\quad =
(2 u_n)^{n-i_n}
\exp \biggl( (n-i_n) \biggl( c_n + \frac{\log(u_n)+2t-\log(\uppi)}{2 u_n}
+ \mathrm{O} \bigl( c_n^2 + \log(u_n)^2
/ u_n^2 \bigr) \biggr) \biggr).
\end{eqnarray*}
Combining (\ref{eqloclevPoisright}), (\ref{expright}) and the last
expression we get
(\ref{alphaextremeright}) in case $n-i_n=\mathrm{o}(\sqrt{u_n})$ and applying
Stirling's formula
(\ref{stirling}) to $(n-i_n)!$ we get (\ref{alphasmallerintermediateright}).
\end{pf*}

\begin{pf*}{Proof of Lemma~\protect\ref{loclevelslargerintermediatecentral}}
Formula (\ref{quotient2}) in Theorem~\ref{zoneofconv} implies
\begin{equation}
\label{BBB}
\alpha_{i_n,n} = \exp\bigl(-x_n^2/2
\bigr) / (\sqrt{2 \uppi}x_n) \bigl[ 1+\mathrm{O} \bigl( 1/x_n^2
+ x_n^3 / \sqrt{n h_{i_n,n} (1-h_{i_n,n})}
\bigr) \bigr].
\end{equation}
First, we have to calculate $x_n$.
From (\ref{hseriescentralintermed}), we get
\[
n h_{i_n,n} (1-h_{i_n,n}) = i_n ( 1- i_n/n
) \bigl[1 + \mathrm{O}\bigl( \sqrt {u_n /\bigl( i_n
(1-i_n/n )\bigr) }\bigr)\bigr]
\]
and hence
\[
\sqrt{ n h_{i_n,n} (1-h_{i_n,n}) } = \sqrt{ i_n ( 1-
i_n/n ) }\bigl[ 1 + \mathrm{O}\bigl( \sqrt{u_n /\bigl(
i_n (1-i_n/n )\bigr) } \bigr) \bigr].
\]
Regarding to (\ref{uniformlyconvergenceofh}), we arrive at
\[
i_n-n h_{i_n,n} = \sqrt{2 i_n (
1-i_n/n) u_n } \bigl[ 1 + \bigl(\log(u_n)+2 t
- \log(\uppi)\bigr)/(4u_n) + \mathrm{O}\bigl( \varepsilon_n(i_n)\bigr)
\bigr],
\]
where $ \varepsilon_n(i_n) = \sqrt{ u_n/( i_n(1-i_n))} + \log(u_n)^2
/ u_n^2 $.
Hence,
\begin{eqnarray*}
x_n & = & \frac{
\sqrt{2 i_n ( 1-i_n/n) u_n }  [
1 + ( \log(u_n)+2 t - \log(\uppi) )/(4u_n) + \mathrm{O}(\varepsilon_n(i_n))
 ] }{
 \sqrt{ i_n ( 1-i_n/n ) } [ 1 + \mathrm{O}(\sqrt{u_n/(i_n(1-i_n/n))}) ] }
\\
& = & \sqrt{2 u_n } \bigl[1 +\bigl(\log(u_n)+2 t-\log(
\uppi)\bigr)/(4u_n) +\mathrm{O}\bigl(\varepsilon _n(i_n)
\bigr) \bigr]
\end{eqnarray*}
and
\[
x_n^2 =2 u_n + \log(u_n)+2 t-
\log(\uppi)+\mathrm{O}\bigl(u_n \varepsilon_n(i_n)\bigr).
\]
This, the fact that
$ 1/x_n = 1/\sqrt{2 u_n} [ 1+ \mathrm{O} ( \log(u_n)/u_n ) ]$ and (\ref{BBB})
lead to
\begin{eqnarray*}
\alpha_{i_n,n} & = & 1/(2\sqrt{\uppi u_n}) \exp \bigl( -
u_n - \log(u_n)/2 - t + \log(\uppi)/2 \bigr)
\\
&&{}\times\bigl[ 1+\mathrm{O}\bigl( \log(u_n)/ u_n +
u_n^{3/2}/\sqrt{i_n(1-i_n/n) }
\bigr) \bigr]
\end{eqnarray*}
and hence (\ref{alphalargerintermediatecentral}) follows.
\end{pf*}

\begin{pf*}{Proof of Lemma~\protect\ref{un3}}
We restrict our attention to $i_n\in B_0\cup B_c$.
The other case can be proved similarly.
Combining (\ref{eqloclevPois2}) and (\ref{stirling}), we get
\begin{equation}
\label{Cn}
\alpha_{i_n,n} = C_n/(2\sqrt{\uppi
u_n}) \bigl[1+\mathrm{o}(1)\bigr]
\end{equation}
with $ C_n \equiv [ (n h_{i_n,n}/i_n) \exp ( 1 - n h_{i_n,n}/
i_n  )  ]^{i_n}$.
It holds
\[
C_n = \exp\bigl( i_n \log\bigl( (n h_{i_n,n} /
i_n) \exp( 1- n h_{i_n,n}/ i_n ) \bigr) \bigr).
\]
From (\ref{hseriescentralintermed}), we get $1-n h_{i_n,n}/i_n =
\mathrm{O}(\sqrt{u_n/i_n})$.
Applying $\log((1-x)\exp(x))=-x^2/2-x^3/3+\mathrm{O}(x^4)$ for $x\in(0,1)$, we
arrive at
\[
C_n = \exp\bigl( i_n \bigl\{ - ( 1- n h_{i_n,n}/
i_n )^2 / 2 - ( 1- n h_{i_n,n}/ i_n
)^3 /3 + \mathrm{O} \bigl( ( u_n/ i_n )^2
\bigr) \bigr\} \bigr).
\]
Lemma~\ref{entwicklungh} leads to
\begin{eqnarray*}
1- n h_{i_n,n}/ i_n &=&  \sqrt{ 2u_n /
i_n } + \frac{\log(u_n)+2t -\log(\uppi)}{2\sqrt{2i_n u_n}} - \frac{u_n}{i_n} + \mathrm{o} (
u_n/ i_n ),
\\
( 1- n h_{i_n,n}/ i_n )^2 &=& \frac{2u_n}{i_n} +
\frac{\log(u_n)+2t
-\log
(\uppi)}{i_n} -\frac{2\sqrt{2} u_n^{3/2}}{i_n^{3/2}} + \mathrm{o} \bigl( u_n^{3/2}/
i_n^{3/2} \bigr)
\end{eqnarray*}
and
\[
( 1- n h_{i_n,n}/ i_n )^3 = 2
\sqrt{2}u_n^{3/2}/ i_n^{3/2} + \mathrm{o} \bigl(
u_n^{3/2} / i_n^{3/2} \bigr).
\]
Then
\begin{eqnarray*}
C_n & = & \exp \biggl( -u_n - \frac{\log(u_n)+2t -\log(\uppi)}{2} +
\frac{\sqrt
{2}u_n^{3/2}}{3\sqrt{i_n}} + \mathrm{o} \biggl( \frac{u_n^{3/2}}{\sqrt{i_n}} \biggr) \biggr)
\\
& = & \exp(-t) \sqrt{\uppi}/\bigl( \log(n) \sqrt{u_n}\bigr) \exp\bigl(
\sqrt {2}u_n^{3/2}/(3\sqrt{i_n}) \bigl(1+ \mathrm{o}(1)
\bigr)\bigr)
\end{eqnarray*}
and hence (\ref{Cn}) yields (\ref{alphaunpower3}) for $i_n$
fulfilling $u_n=\mathrm{o}(i_n)$ and $i_n=\mathrm{O}(u_n^3)$.
\end{pf*}

\begin{pf*}{Proof of Lemma~\protect\ref{loclevelnurun}}
Formulas (\ref{inCun}) and (\ref{ninCun}) provide that in order
to find $\alpha_{i_n,n}$
we have to calculate $\mathbb{P}(Y_n=i_n)$ and $\mathbb{P}(\tilde
{Y}_n=n-i_n)$.

We start with the case $c_n\equiv u_n/i_n\to c>0$.
Noting that $(1-x)^{k}=\exp(-k x+\mathrm{O}(k x^2))$ for $x\in( 0,1)$ and
$k\in
\mathbb{N}$, formula (\ref{hseriesintermediate}) implies
\[
\biggl( \frac{n h_{i_n,n}}{i_n} \biggr)^{i_n} = \bigl(\delta(c_n)
\bigr)^{i_n} \exp \biggl( - \frac{\log(u_n)+2 t-\log(\uppi)}{2 \sqrt{c_n^2+2c_n}} + \mathrm{O} \biggl(
\frac
{\log(u_n)^2}{u_n} \biggr) \biggr).
\]
Applying the Stirling formula (\ref{stirling}) and (\ref
{hseriesintermediate}), we arrive at
\begin{eqnarray*}
\mathbb{P}(Y_n=i_n) & = & ( n h_{i_n,n} /
i_n )^{i_n} /\sqrt{2\uppi i_n}
\exp(i_n-n h_{i_n,n})
\\
& = & \bigl[ \delta(c_n) \exp\bigl(1-\delta(c_n)\bigr)
\bigr]^{i_n} /\sqrt{2\uppi i_n} \bigl( \sqrt{\uppi}\exp(-t) /
\sqrt{u_n} \bigr)^{({1-\delta(c_n)})/{\sqrt{c_n^2+2c_n}}}
\\
&& {}\times\bigl( 1+\mathrm{O} \bigl( \log(u_n)^2 /
u_n \bigr) \bigr).
\end{eqnarray*}
Therefore, (\ref{inCun}) implies (\ref{alphanurunleft}).

Now let $i_n$ be such that $c_n\equiv(n-i_n)/u_n\to c>0$.
Similarly as above, (\ref{hseriesextremeright}) implies
\begin{eqnarray*}
\biggl( \frac{n-n h_{i_n,n}}{n-i_n} \biggr)^{n-i_n} & = & \bigl[(1+c_n+
\sqrt{1+2 c_n})/c_n \bigr]^{n-i_n}
\\
&& {}\times\exp\bigl(
c_n/\sqrt{1+2 c_n}\bigl(\log(u_n)+2 t -\log(\uppi)\bigr)/2 + \mathrm{O} \bigl(
\log(u_n)^2/u_n \bigr) \bigr).
\end{eqnarray*}
Then
\begin{eqnarray*}
\mathbb{P}(\tilde{Y}_n=n-i_n) & = & \biggl[
\frac{n-n h_{i_n,n}}{n-i_n} \biggr]^{n-i_n} \Big/ \sqrt{2 \uppi(n-i_n)}
\exp(n h_{i_n,n}-i_n)
\\
& = & \bigl[ (1+c_n+\sqrt{1+2 c_n})/c_n
\bigr]^{n-i_n} /\sqrt{2 \uppi(n-i_n)}
\\
&& {}\times\bigl( \sqrt{\uppi}\exp(-t)/ \sqrt{u_n}
\bigr)^{1+{1}/{\sqrt{1+2 c_n}}} \log(n)^{-(1+\sqrt{1+2 c_n})}
\\
&& {}\times\bigl( 1 + \mathrm{O} \bigl( \log(u_n)^2 /
u_n \bigr) \bigr)
\end{eqnarray*}
and (\ref{ninCun}) lead to (\ref{alphanurunright}).
\end{pf*}

\section{Proofs of Section~\texorpdfstring{\protect\ref{secquotient}}{5}}\label{AppC}

\begin{pf*}{Proof of Theorems \protect\ref{locallevelscomparisons1}}
Formula (\ref{aim1}) for the case (i) immediately\vspace*{1pt} follows from
Lemma~\ref{loclevelslargerintermediatecentral}.
Here we prove (\ref{aim1}) for (ii)--(v), (\ref{aim01}) for (vi)--(ix), (\ref{aim0}) for
$(i_n^{(1)},i_n^{(2)})\in B_0 \times B_0$ such that (iv), (viii) are not fulfilled and
(\ref{aim0}) for $(i_n^{(1)},i_n^{(2)})\in\bar{B}_0 \times\bar{B}_0$
such that\vspace*{1pt} (v), (ix)
are not fulfilled. Lemma~\ref{lemun3} shows (\ref{aim01}) for (x),
(xi),\vspace*{1pt} (\ref{aim0}) for
$(i_n^{(1)},i_n^{(2)})\in A_c \times A_c$, $c>0$, such that (x) is
not fulfilled and (\ref{aim0})
for $(i_n^{(1)},i_n^{(2)})\in\bar{A}_c \times\bar{A}_c$, $c>0$ such
that (xi) is not fulfilled.
The remaining cases for (\ref{aim0}) are proved in Lemmas \ref{lemi1-in}, \ref{compare1},
\ref{compare2}, \ref{compare3} and \ref{compare4}.

For $(i_n^{(1)},i_n^{(2)})\in(B_0 \cup B_c)\times(B_0 \cup B_c)$
Lemma~\ref{un3} yields
\[
\alpha_{i_n^{(2)},n} / \alpha_{i_n^{(1)},n} = \exp\Bigl( - \sqrt{2}
u_n^{3/2} \big/\Bigl( 3 \sqrt{i_n^{(1)}}
\Bigr) \Bigl( 1 - \sqrt{ i_n^{(1)} /
i_n^{(2)}} \Bigr) \bigl[1+\mathrm{o}(1)\bigr] \Bigr).
\]
This implies (\ref{aim1}) for (ii) and (\ref{aim01}) for
$(i_n^{(1)},i_n^{(2)})\in B_{c_1} \times B_{c_2}$
with $c_1<c_2$, that is, (\ref{aim01}) for a partial case of (vi).
Moreover, Lemmas \ref{loclevelslargerintermediatecentral} and \ref{un3} immediately yield the remaining cases of~(vi).

For $i_n^{(1)}\in B_0$ define $b_n \equiv u_n^{3/2} /\sqrt{i_n^{(1)}}
(1-\sqrt{i_n^{(1)}/i_n^{(2)}}) $.
Clearly, we get (\ref{aim1}) if $b_n\to0$, (\ref{aim01}) if $b_n \to
c>0$ and (\ref{aim0}) in case $b_n\to\infty$.
Note that
\[
i_n^{(1)} / i_n^{(2)} =
1-b_n \sqrt{ i_n^{(1)}
}/u_n^{3/2} \Bigl( 2 + b_n \sqrt{
i_n^{(1)}} / u_n^{3/2} \Bigr).
\]
If $b_n\to0$ as $n\to\infty$, that is, $i_n^{(1)}/i_n^{(2)} = 1 +
\mathrm{o}(\sqrt{i_n^{(1)}}/u_n^{3/2})$, then we get~(\ref{aim1}) for (iv).
We get~(\ref{aim01}) for (viii), when $b_n\to b$ for some $b>0$
and (\ref{aim0}) in case $b_n\to\infty$.

Finally, (iii), (v), (vii) and (ix) can be proved in a
similar way.
\end{pf*}

\begin{lemma}\label{lemi1-in}
Let\vspace*{1pt} $\{ i_n^{(1)} \}_{n\in\mathbb{N}}$ and $\{ i_n^{(2)} \}_{n\in
\mathbb{N}}$ be
such that $i_n^{(1)}<i_n^{(2)}$,
$n\in\mathbb{N}$, and either $(i_n^{(1)},i_n^{(2)}) \in A_0\times
A_0$ or
$(i_n^{(1)},i_n^{(2)}) \in\bar{A}_0\times\bar{A}_0$.
Then (\ref{aim0}) is fulfilled.
\end{lemma}

\begin{pf}
First, let $i_n^{(j)}\in A_0$, $j=1,2$.
If $i_n^{(2)}\equiv i_2\in\mathbb{N}$ for all $n\in\mathbb{N}$ and
$n$ is large
enough, representation (\ref{alphaextremeleft})
immediately yields $\alpha_{i_2,n} / \alpha_{i_1,n} = \mathrm{O} ( (2
u_n)^{i_1-i_2}  )$, and hence (\ref{aim0}) follows.
Furthermore, let $i_n^{(2)}\to\infty$ as $n\to\infty$.
Since $1/\sqrt{2 \uppi i}$ and $\exp(-i v_n)$ in representation
(\ref{alphasmallerintermediateleft})
decrease as $i$ increases for a fixed larger $n$, in order to prove
(\ref{aim0}) it suffices to show that
\[
B_n \equiv \log\bigl( \bigl( \gamma i_n^{(2)}/
u_n \bigr)^{i_n^{(2)}} / \bigl( \gamma i_n^{(1)}/
u_n \bigr)^{i_n^{(1)}} \bigr)
\]
converges to $-\infty$ as $n\to\infty$.
Setting $x_n\equiv i_n^{(2)}-i_n^{(1)}$, we obtain
\[
B_n = i_n^{(1)} \bigl[ \bigl(- x_n/i_n^{(1)}\bigr)
\log( u_n / \gamma) + \bigl( 1 + x_n /
i_n^{(1)} \bigr) \log \bigl(i_n^{(2)}
\bigr)- \log\bigl(i_n^{(1)}\bigr) \bigr].
\]
If $d_n\equiv x_n/i_n^{(1)} \to d $ for a $d>0$ or $d=\infty$, we get
$B_n = - i_n^{(1)} d_n \log(u_n/\gamma)(1 + \mathrm{o}(1)) $, that is, $B_n\to
-\infty$ as $n\to\infty$.
Hence, (\ref{aim0}) is fulfilled.

For $i_n^{(j)}$, $j=1,2$, such that $x_n/i_n^{(1)} =\mathrm{o}(1) $ we get
\[
B_n = - x_n \log( u_n/\gamma) +
x_n \log\bigl(i_n^{(2)}\bigr)
+i_n^{(1)} \log\bigl( 1+ x_n/
i_n^{(1)} \bigr).
\]
Applying $\log(1+x)=x+\mathrm{O}(x^2)$ for $x\in(0,1)$ and the fact that
$i_n^{(2)} = \mathrm{o}(u_n)$, we obtain
$ B_n = - x_n \log( u_n/ \gamma)(1+\mathrm{o}(1)) $, and hence (\ref
{aim0}) follows.

Now, let $i_n^{(j)} \in\bar{A}_0$, $j=1,2$.
For $i_n^{(1)}< i_n^{(2)}$ such that $n-i_n^{(1)}$ is fixed, formula
(\ref{alphaextremeright}) in
Lemma~\ref{loclevelsextremeintermediateright} immediately leads to
the assertion.
For the case $n-i_n^{(1)}\to\infty$ as $n\to\infty$, due to~(\ref{alphasmallerintermediateright})
it suffices to consider
\[
D_n\equiv \bigl( u_n/\bigl(\gamma(n-i_n)
\bigr) \bigr)^{n-i_n} / \sqrt{n-i_n} \exp\bigl((n-i_n)
w_n\bigr).
\]
Since
\[
D_n = \exp \bigl( (n-i_n) \bigl[ - \log(n-i_n)
/ \bigl(2(n-i_n)\bigr) + \log \bigl( u_n/\bigl(
\gamma(n-i_n)\bigr) \bigr) +w_n \bigr] \bigr),
\]
$\log(u_n/(\gamma(n-i_n)))\to\infty$ as $n\to\infty$, $\log(x)/x<1$
for $x\geq1$ and $w_n=\mathrm{o}(1)$,
we arrive at $ D_n = \exp ( (n-i_n) \log ( u_n/(\gamma(n-i_n))
 ) [1+\mathrm{o}(1)]  ) $.
Thus, it suffices to show that
\[
\log\bigl( \bigl( u_n/\bigl(\gamma\bigl(n-i_n^{(2)}
\bigr)\bigr) \bigr)^{n-i_n^{(2)}} / \bigl( u_n/\bigl( \gamma
\bigl(n-i_n^{(1)}\bigr)\bigr) \bigr)^{n-i_n^{(1)}} \bigr)
\]
converges to $-\infty$ for $n\to\infty$.
This can be proved similarly as before.
\end{pf}

\begin{lemma}\label{lemun3}
Let\vspace*{1pt} $\{ i_n^{(1)} \}_{n\in\mathbb{N}}$ and $\{ i_n^{(2)} \}_{n\in
\mathbb{N}}$ be
such that
$i_n^{(1)}<i_n^{(2)}$ for $n\in\mathbb{N}$. We suppose that either
$\lim_{n\to\infty}
u_n/i_n^{(j)} = c_j$, $j=1,2$,
or $\lim_{n\to\infty}(n-i_n^{(j)})=c_j $, $j=1,2$, for arbitrary but
fixed $c_j>0$.
Moreover, let $m_n\equiv i_n^{(2)}-i_n^{(1)}$, $n\in\mathbb{N}$.
If $m_n = m$ for some fixed $m\in\mathbb{N}$ and all $n\in\mathbb
{N}$, (\ref{aim01}) is
fulfilled and if $m_n \to\infty$ as $n\to\infty$, (\ref{aim0}) is
fulfilled.
\end{lemma}

\begin{pf}
First, let $c_{j,n}\equiv u_n/i_n^{(j)} \to c_j>0$, $j=1,2$.
For $i_n$ such that $c_n\equiv u_n/i_n\to c>0$ formula~(\ref
{alphanurunleft}) in Lemma~\ref{loclevelnurun} implies
\[
\alpha_{i_n,n} = \exp \bigl( f_1(c_n)
u_n +f_2(c_n)\log(u_n)
+f_3(c_n) + \mathrm{o}(1) \bigr),
\]
where
\begin{eqnarray*}
f_1(c) &=& (1/c) \log \bigl( \delta(c) / \exp\bigl(\delta(c)-1\bigr)
\bigr),\qquad
f_2(c) = -1+c/\bigl(2 \sqrt{c^2+2c}\bigr),
\\
f_3(c) &=& \log \bigl( \sqrt{c}/\bigl(\bigl(1-\delta(c)\bigr)\sqrt{2
\uppi}\bigr) \bigr) + \bigl(1-\delta(c)\bigr)/\sqrt{c^2+2c} \log\bigl(
\sqrt{\uppi} \exp(-t)\bigr).
\end{eqnarray*}
It follows
%
\begin{eqnarray*}
\alpha_{i_n^{(2)},n} / \alpha_{i_n^{(1)},n} & = & \exp\bigl(
\bigl(f_1(c_{2,n})-f_1(c_{1,n})\bigr)
u_n + \bigl(f_2(c_{2,n})-f_2(c_{1,n})
\bigr)\log(u_n)
\\
&&\hspace*{18pt}{}+ f_3(c_{2,n})-f_3(c_{1,n})
+ \mathrm{o}(1)\bigr).
\end{eqnarray*}
%
Since $f_1(\cdot)$ is strictly increasing, (\ref{aim0}) is fulfilled
in case $c_1>c_2$.
Now, let $c_1=c_2$.
Hence,
\[
\frac{ \alpha_{i_n^{(2)},n} }{ \alpha_{i_n^{(1)},n} } = \exp \bigl( \bigl(f_1(c_{2,n})-f_1(c_{1,n})
\bigr) u_n + \bigl(f_2(c_{2,n})-f_2(c_{1,n})
\bigr) \log(u_n) +\mathrm{o}(1) \bigr).
\]
Setting $x_n\equiv i_n^{(2)}-i_n^{(1)}$ and noting that $x_n=\mathrm{o}(u_n)$,
we get
\begin{eqnarray*}
c_{1,n}&=& c_{2,n}/\bigl(1-x_n/i_n^{(2)}
\bigr) = c_{2,n} \bigl[1 + c_{2,n} x_n /
u_n + \mathrm{O} \bigl( x_n^2 / u_n^2
\bigr) \bigr],
\\
f_1(c_{1,n}) &=& f_1(c_{2,n}) +
c_{2,n} \Bigl[ c_{2,n}/\sqrt {c_{2,n}^2+2c_{2,n}}-1
- f_1(c_{2,n}) \Bigr] x_n/u_n +\mathrm{O} \bigl( x_n^2/u_n^2 \bigr)
\end{eqnarray*}
and $f_2(c_{1,n}) = f_2(c_{2,n}) + \mathrm{O} ( x_n/u_n ) $.
Therefore,
\[
\alpha_{i_n^{(2)},n} / \alpha_{i_n^{(1)},n} = \exp\Bigl( x_n
c_{2,n} \Bigl[ 1+f_1(c_{2,n})- c_{2,n}\big/
\sqrt{c_{2,n}^2+2 c_{2,n}} +\mathrm{o}(1) \Bigr]
\Bigr).
\]
Since $1+f_1(c)-c/\sqrt{c^2+2c} < 0$ for $c\in(0,\infty)$, we get
(\ref{aim01}) if $x_n=x$ for
some $x\in\mathbb{N}$ and~(\ref{aim0}) if $x_n\to\infty$ as $n\to
\infty$.

The case $c_{j,n} \equiv(n-i_n^{(j)})/u_n\to c_j$, $j=1,2$, can be
proved similarly.
\end{pf}

\begin{lemma}\label{compare1}
Let $\{i_n^{(1)}\}_{n\in\mathbb{N}}$ and $\{i_n^{(2)}\}_{n\in\mathbb
{N}}$ be such
that $i_n^{(1)}=\mathrm{o}(u_n)$ and
$c_n \equiv u_n/i_n^{(2)}\to c$, $n\to\infty$, for some $c>0$.
Then (\ref{aim0}) is fulfilled.
\end{lemma}

\begin{pf}
Formulas (\ref{alphasmallerintermediateleft}) and (\ref
{alphanurunleft}) yield
\begin{eqnarray*}
\log \biggl( \frac{\alpha_{i_n^{(2)},n}}{\alpha_{i_n^{(1)},n}} \biggr) & = & u_n \biggl[ - \biggl(
\frac{1}{2}+\frac{1-\delta(c_n)}{2\sqrt{c_n^2+2 c_n}} \biggr) \frac{\log(u_n)}{u_n} +
\frac{\log(i_n^{(1)})}{2u_n} +\frac{i_n^{(1)}}{u_n} v_n
\\
&& \hspace*{7pt}\quad{}+\frac{1}{c_n} \log \biggl( \frac{\delta(c_n)}{\exp(\delta(c_n)-1)} \biggr) -
\frac{i_n^{(1)}}{u_n} \log \biggl( \gamma\frac{i_n^{(1)}}{u_n} \biggr) + \mathrm{o}(1) \biggr].
\end{eqnarray*}
Since $i_n^{(1)}/u_n=\mathrm{o}(1)$ and $\lim_{x\to0} x \log(x) =0$ and $\log
(\delta(c)/\exp(\delta(c)-1))<0$
for all $c\in(0,\infty)$, we obtain $\log ( \alpha
_{i_n^{(2)},n} /
\alpha_{i_n^{(1)},n}  )
\to-\infty$ as $n\to\infty$ and hence (\ref{aim0}) follows.
\end{pf}

\begin{lemma}\label{compare2}
Let\hspace*{1pt} $\{i_n^{(1)}\}_{n\in\mathbb{N}}$ and $\{i_n^{(2)}\}_{n\in\mathbb
{N}}$ be such
that $c_n \equiv u_n/i_n^{(1)}\to c$,
$n\to\infty$, for some $c>0$, $i_n^{(2)}=\mathrm{O}(u_n^3)$ and $u_n=\mathrm{o}(i_n^{(2)})$.
Then (\ref{aim0}) is fulfilled.
\end{lemma}

\begin{pf}
Formulas (\ref{alphaunpower3}) and (\ref{alphanurunleft}) imply
\[
\log ( \alpha_{i_n^{(2)},n} / \alpha_{i_n^{(1)},n} ) = - \bigl[
1+(1/c_n) \log \bigl( \delta(c_n)/ \exp\bigl(
\delta(c_n)-1\bigr) \bigr) \bigr] u_n +\mathrm{o}(u_n).
\]
Since $c_n\to c>0$ as $n\to\infty$ and $\log ( \delta(c)/\exp
(\delta
(c)-1)  )/c >-1$
for all $c\in(0,\infty)$, we immediately obtain (\ref{aim0}).
\end{pf}

\begin{lemma}\label{compare3}
Let $\{i_n^{(1)}\}_{n\in\mathbb{N}}$ and $\{i_n^{(2)}\}_{n\in\mathbb
{N}}$ be such
that $n-i_n^{(1)}=\mathrm{O}(u_n^3)$,
$u_n=\mathrm{o}(n-i_n^{(1)})$ and $c_n \equiv(n-i_n^{(2)})/u_n\to c$, $n\to
\infty$, for some $c>0$.
Then (\ref{aim0}) is fulfilled.
\end{lemma}

\begin{pf}
Formulas (\ref{alphaunpower3}) and (\ref{alphanurunright}) lead to
\[
\log \biggl( \frac{\alpha_{i_n^{(2)},n}}{\alpha_{i_n^{(1)},n}} \biggr) = - \biggl[ \sqrt{1+2 c_n}
- c_n \log \biggl( \frac{1+c_n+\sqrt{1+2
c_n}}{c_n} \biggr) \biggr]
u_n + \mathrm{o}(u_n).
\]
Noting that $\sqrt{1+2 c} - c \log((1+c+\sqrt{1+2 c})/c)>0$ for all
$c\in(0,\infty)$, we get (\ref{aim0}).
\end{pf}

\begin{lemma}\label{compare4}
Let $\{i_n^{(1)}\}_{n\in\mathbb{N}}$ and $\{i_n^{(2)}\}_{n\in\mathbb
{N}}$ satisfy
$c_n \equiv(n-i_n^{(1)})/u_n\to c$,
$n\to\infty$, for some $c>0$ and $n-i_n^{(2)}=\mathrm{o}(u_n)$. Then (\ref
{aim0}) is fulfilled.
\end{lemma}

\begin{pf}
Formulas (\ref{alphasmallerintermediateright}) and (\ref
{alphanurunright}) lead to
\[
\log \biggl( \frac{\alpha_{i_n^{(2)},n}}{\alpha_{i_n^{(1)},n}} \biggr) = \biggl[ - 1 + \sqrt{1+2
c_n} - c_n \log \biggl( \frac{1+c_n+\sqrt{1+2
c_n}}{c_n} \biggr)
\biggr] u_n +\mathrm{o}(u_n).
\]
Since $-1+\sqrt{1+2 c} - c \log((1+c+\sqrt{1+2 c})/c)< 0$ for all
$c\in
(0,\infty)$, (\ref{aim0}) is fulfilled.
\end{pf}
\end{appendix}

\section*{Acknowledgements}
The authors sincerely thank the
referees
for numerous helpful and constructive comments and suggestions and
several hints for additional references.
Special thanks are due to the Editor E.~Moulines for his patience in
handling the manuscript.
This work was supported by the Ministry of Science and Research of the
State of North Rhine-Westphalia (MIWF NRW) and the German Federal
Ministry of Health (BMG).






\printhistory
\end{document}